\documentclass[review,onefignum,onetabnum]{siamart190516}

\usepackage{lipsum}
\usepackage{graphicx}
\usepackage{graphicx}
\usepackage{epstopdf}
\usepackage{algorithmic}
\usepackage{amssymb}
\usepackage[export]{adjustbox}
\ifpdf
  \DeclareGraphicsExtensions{.eps,.pdf,.png,.jpg}
\else
  \DeclareGraphicsExtensions{.eps}
\fi

\usepackage{caption}

 \usepackage{mathtools}

\newsiamremark{remark}{Remark}
\newsiamremark{hypothesis}{Hypothesis}
\crefname{hypothesis}{Hypothesis}{Hypotheses}
\newsiamthm{claim}{Claim}

\usepackage{amsfonts}
\DeclareMathOperator*{\argmin}{\mbox{arg}\min}

\makeatletter
\newcommand{\setword}[2]{%
  \phantomsection
#1\def\@currentlabel{\unexpanded{#1}}\label{#2}%
}
\makeatother

\usepackage{lipsum}
\usepackage{amsfonts}
\usepackage{graphicx}
\usepackage{epstopdf}
\ifpdf
  \DeclareGraphicsExtensions{.eps,.pdf,.png,.jpg}
\else
  \DeclareGraphicsExtensions{.eps}
\fi

\headers{Dewetting}{J. Guo, S. Esedo\=glu}

\title{Median Filters for\\ Anisotropic Wetting / Dewetting Problems\thanks{Submitted to the editors June 19, 2024.\funding{Jiajia Guo and Selim Esedo\=glu acknowledge support from the NSF-DMS 2012015 and NSF-DMS 2410272.}}}

\author{Jiajia Guo \thanks{Corresponding author. 530 Church St, Ann Arbor, MI 48109. {\bf jiajiag@umich.edu}}\\ University of Michigan \and Selim Esedo\=glu\\ University of Michigan}

\usepackage{amsopn}

\ifpdf
\hypersetup{
  pdftitle={Median Filters for\\ Anisotropic Wetting / Dewetting Problems},
  pdfauthor={Jiajia Guo \thanks{Corresponding author. 530 Church St, Ann Arbor, MI 48109. {\bf jiajiag@umich.edu}}\\ University of Michigan \and Selim Esedo\=glu\\ University of Michigan}
}
\fi

\begin{document}
\nolinenumbers
\maketitle

\begin{abstract}
 \noindent We present new level set methods for multiphase, anisotropic (weighted) motion by mean curvature of networks, focusing on wetting-dewetting problems where one out of three phases is stationary -- a good testbed for checking whether complicated junction conditions are correctly enforced.
The new schemes are vectorial median filters: The level set values at the next time step are determined by a sorting procedure performed on the most recent level set values.
Detailed numerical convergence studies are presented, showing that the correct angle conditions at triple junctions (which include torque terms due to anisotropy) are indeed indirectly and automatically attained.
Other standard benefits of level set methods, such as subgrid accuracy on uniform grids via interpolation and seamless treatment of topological changes, remain intact.
\end{abstract}



\section{Introduction and Background}
\let\thefootnote\relax\footnotetext{{\it 2020 Mathematics Subject Classification.} Primary: 65M12; Secondary 35K93.}
Multiphase, anisotropic (weighted) motion by mean curvature of networks of interfaces (curves in $\mathbb{R}^2$, surfaces in $\mathbb{R}^3$) is a challenging dynamics to simulate that arises in many important applied contexts, such as coarsening of microstructure in polycrystalline materials \cite{kinderlehrer1} that include most metals and ceramics.
The evolution is variational: It can be viewed, at least formally, as gradient descent for a cost function, or energy, defined on partitions of a domain $D \in \mathbb{R}^d$ into disjoint sets.
In the generality of interest in this paper, the energy is given by:
\begin{equation}
\label{eq:multiphase}
\mathcal{E}(\Sigma_1,\Sigma_2,\ldots,\Sigma_N) = \sum_{i\not=j} \int_{\Gamma_{i,j}} \sigma_{i,j}(n(x)) \, dS(x)
\end{equation}
where the sets $\Sigma_1,\ldots,\Sigma_N$ partition a domain $D \subset \mathbb{R}^d$ so that
\begin{equation}
\label{eq:partition}
\Gamma_{i,j} := (\partial\Sigma_i) \cap (\partial\Sigma_j) = \Sigma_i \cap \Sigma_j \mbox{ for $i\not=j$, and } \bigcup_{i=1}^N \Sigma_i = D
\end{equation}
and $n(x)$ denotes a unit normal to the interface $\Gamma_{i,j}$ for a point $x$ taken on it.
The continuous functions $\sigma_{i,j} : \mathbb{S}^{d-1} \to \mathbb{R}^+\setminus\{0\}$ are the {\em surface tensions} associated with the interfaces, and are required to be centrally symmetric (even) and have a {\em convex} one-homogeneous extension to $\mathbb{R}^d$:
\begin{equation}
\sigma_{i,j}(x) = |x| \sigma_{i,j} \left( \frac{x}{|x|}\right) \mbox{ for } x\not=0.
\end{equation}
Thus, each surface tension $\sigma$ is a norm on $\mathbb{R}^d$.
Its unit ball $B_\sigma = \{ x \, : \, \sigma(\xi) \leq 1 \}$ is often called the {\em Frank diagram}.
The unit ball $W_\sigma$ of its dual $\sigma^*$, defined as
\begin{equation}
\label{eq:convexsigma}
\sigma^*(\xi) = \sup_{x\in\mathbb{R}^d}  x \cdot \xi
\end{equation}
is known as the {\em Wulff shape} associated with the anisotropy $\sigma$.
A fundamental fact is that $W_\sigma$ solves the isoperimetric problem with anisotropy $\sigma$:
\begin{equation}
W_\sigma \in \mbox{arg-}\min_{\substack{\Sigma \subset \mathbb{R}^d\\ |\Sigma| = C}} \int_{\partial\Sigma} \sigma(n(x)) \, dS(x)
\end{equation}
for an appropriate scalar $C>0$, where we write $|\Sigma|$ to denote the volume of $\Sigma$.
We will be content with $B_\sigma$, and consequently $W_\sigma$, that are strongly convex and have smooth boundaries, thus staying away from the particularly challenging crystalline case (except by approximation).
Well-posedness of model (\ref{eq:multiphase}) also requires the pointwise triangle inequality
\begin{equation}
\label{eq:triangle}
\sigma_{i,j}(n) + \sigma_{j,k}(n) \geq \sigma_{i,k}(n)
\end{equation}
for every $n\in\mathbb{S}^{d-1}$ and distinct $i$, $j$, and $k$.
\textcolor{black}{If this condition is violated, variational model (\ref{eq:multiphase}) fails to be lower semicontinuous, which manifests itself in phase $j$ wetting the interface between phases $i$ and $k$, leading to an {\em effective} surface tension lower than $\sigma_{i,k}$}; see e.g. \cite{morgan1997}.
It is known that these conditions are necessary, but in general not sufficient (except for $N=2$) for well-posedness of the multiphase model (\ref{eq:multiphase}); see e.g. \cite{ambrosio_braides_I, ambrosio_braides_II}.

Under gradient flow for energy (\ref{eq:multiphase}), the motion of the interface $\Gamma_{i,j} \subset \mathbb{R}^3$ at a point $x\in\Gamma_{i,j}$ away from {\em junctions} (along which three or more interfaces meet) is described by {\em weighted motion by mean curvature}, which specifies the normal speed
\begin{multline}
\label{eq:wk}
v_\perp(x) = \Big( \partial_{s_1}^2 \sigma_{i,j} \big( n_{i,j}(x) \big) + \sigma_{i,j} \big( n_{i,j}(x) \big) \Big) \kappa_1(x)\\
 + \Big( \partial_{s_2}^2 \sigma_{i,j} \big( n_{i,j}(x) \big) + \sigma_{i,j} \big( n_{i,j}(x) \big) \Big) \kappa_2(x) 
\end{multline}
where $\kappa_1$ and $\kappa_2$ are the two principal curvatures, and $\partial_{s_\ell}$ denotes differentiation along the great circle on $\mathbb{S}^2$ that passes through the unit normal $n_{i,j}(x)$ at $x$ and has as its tangent the $\ell$-th principal curvature direction.
In the isotropic setting, this simplifies to
\begin{equation}
\label{eq:wksimple}
v_\perp(x) = \sigma_{i,j} \big( \kappa_1(x) + \kappa_2(x) \big).
\end{equation}
In addition, there are what are known as Herring angle conditions \cite{herring} that hold along free boundaries, called triple junctions, formed by the intersection of three interfaces;.
\begin{figure}[h]
\begin{center}
\includegraphics[scale=0.9]{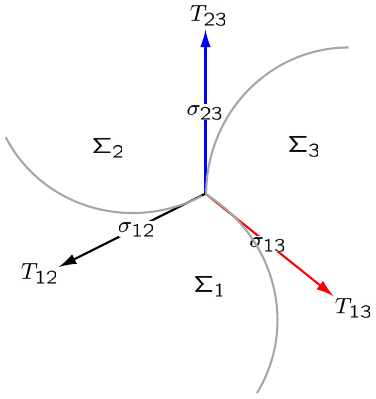}
\caption{\footnotesize A triple junction in two dimensions}
\label{fig:illustration1}
\end{center}
\end{figure}
In two dimensions, these read
\begin{multline}
\label{eq:herringaniso}
\sigma_{1,2}(T^\perp_{1,2})T_{1,2} + \sigma_{1,3}(T^\perp_{1,3})T_{1,3} + \sigma_{2,3}(T^\perp_{2,3})T_{2,3}\\
+\sigma'_{1,2}(T^\perp_{1,2})T^\perp_{1,2} + \sigma'_{1,3}(T^\perp_{1,3})T^\perp_{1,3} + \sigma'_{2,3}(T^\perp_{2,3})T^\perp_{2,3} = 0.
\end{multline}
where $T_{i,j}$ is the unit tangent to interface $\Gamma_{i,j}$ at the triple junction, pointing away from the junction; see Figure \ref{fig:illustration1}.
In the case of the simpler isotropic model where $\sigma_{i,j}$ in (\ref{eq:multiphase}) are {\em constants}, it simplifies to 
$$\sigma_{1,2} T_{1,2} + \sigma_{1,3} T_{1,3} + \sigma_{2,3} T_{2,3} = 0 $$
$$ \iff \frac{\sin\theta_1}{\sigma_{2,3}} = \frac{\sin\theta_2}{\sigma_{1,3}} = \frac{\sin\theta_3}{\sigma_{1,2}} $$
so that the angles formed at the junction are determined in terms of the surface tensions.

Two essential difficulties in developing algorithms for model (\ref{eq:wk}) \& (\ref{eq:herringaniso}) are {\bf 1.} topological changes in the evolving network of \textcolor{black}{\textbf{boundaries}} \textcolor{black}{$\cup\Gamma_{i,j}$}, and {\bf 2.} imposition of the Herring angle conditions (\ref{eq:herringaniso}) along junctions which are themselves part of the unknown.
The first difficulty makes numerical methods that represent interfaces implicitly, such as the level set method \cite{osher_sethian}, particularly appealing.
However, boundary conditions such as (\ref{eq:herringaniso}) are not imposed explicitly in these methods; instead, it is necessary to ensure they are induced indirectly without compromising the simplicity and elegance that are a defining feature of these numerical schemes.
There is as of yet no level set method that can handle the full generality of model (\ref{eq:wk}) \& (\ref{eq:herringaniso})  in a robust, reliable way, guaranteeing in particular the all important junction conditions.

A recent approach \cite{esedoglu_guo_li} to numerical approximation of multiphase motion by mean curvature in the isotropic case (where surface tensions $\sigma_{i,j}$ are constants) utilizes vectorial versions of the median filter.
It is an extension to the multiphase setting of the median filter scheme that was proposed for two-phase motion by mean curvature in \cite{oberman2004} as a discretization of the standard level set PDE for this geometric motion.
In this paper, we take a step towards further extension of median filter schemes to the {\em anisotropic} multiphase motion by mean curvature context of model (\ref{eq:multiphase}).
In particular, we focus on the relatively simpler wetting / dewetting problems \cite{xu_wang} in two dimensions that feature three phases (liquid, solid, vapor), one of which (the solid phase) is stationary, making them quite close to two-phase problems, while still featuring triple junction conditions.
As the class of numerical schemes we study induce the angle conditions at junctions indirectly, convergence to the correct evolution is already an interesting question in this simpler setting; see e.g. \cite{svadlenka2022} that numerically investigates convergence of several existing approaches in this case.


\section{The \textcolor{black}{Vapor-Liquid-Substrate} Model in 2D}
\begin{figure}
\begin{center}
\includegraphics[scale=0.6]{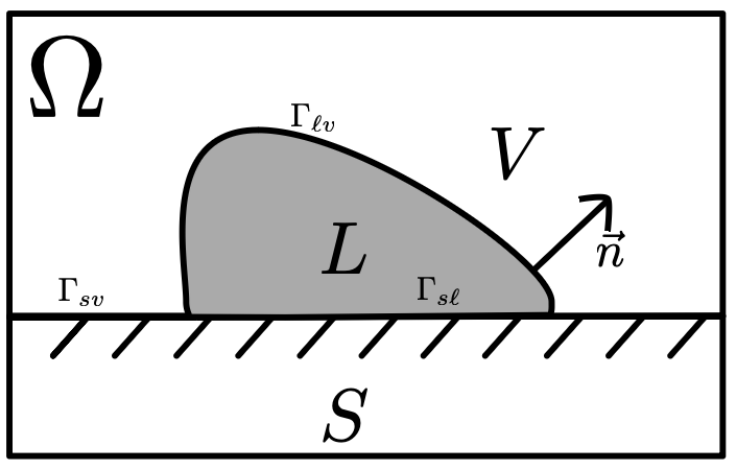}
\caption{\footnotesize Diagram of the dewetting/wetting setting. \textcolor{black}{The substrate $S$ need not be flat.}}
\label{fig:diagram1}
\end{center}
\end{figure}

\textcolor{black}{
We consider models of a liquid drop $L$ spreading on a solid surface $S$ under surface tension effects, in the presence of a third, vapor phase $V$ that surrounds them in a bounded domain $\Omega \in\mathbb{R}^2$, where $V, L, S$ are essentially disjoint subsets of $\Omega$, with $|V|+|L|+|S|=|\Omega|$.
Here, for a set $\Sigma$ we denote its \textcolor{black}{area} by $|\Sigma|$. \textcolor{black}{See Figure \ref{fig:diagram1} for an illustration}. 
Throughout the paper, we impose periodic boundary conditions for convenience, so that $\Omega$ is in fact a torus.}
The interfaces separating the three phases $V$, $L$, and $S$ will be denoted by \textcolor{black}{$\Gamma_{VL}$, $\Gamma_{LS}$ and $\Gamma_{VS}$}. 
The solid region described by the set $S$ is stationary; only the sets $L$ and $V$ evolve, partitioning at any time the complement $S^c \cap \Omega$ of $S$.
Hence, $L$ by itself is sufficient to determine the partition of $\Omega$. 
Area of the liquid phase droplet $L$ is preserved during the evolution; we write $|L|=A$ to denote this fixed area.

\textcolor{black}{
A similar set up, where three disjoint phases that preserve their respective areas partition a domain, and where one of the phases is stationary, arises in solid-state dewetting problems \cite{thompson,jiang_wang,bao_jiang,bao_jiang2}; in this case, the surface tension of the solid ``drop'' and the vapor phase may well be anisotropic.
In these applications, the motion of the interface between the drop (or particle) $L$ and the vapor phase $V$ is often modeled by surface diffusion -- a high order geometric flow that is challenging to approximate by numerical techniques that represent interfaces implicitly and thus allow seamless handling of topological changes.
Stationary states of these models can also be sought via $L^2$ steepest descent, which leads to the simpler, second order geometric motion known as (weighted) motion by mean curvature discussed in the previous section.
Motion by mean curvature is also the appropriate dynamics for mesoscale models of grain boundary motion in polycrystalline materials; as already mentioned, in that context the present set up, where one of the interfaces is stationary (i.e. of very low mobility) constitutes a simplified testbed for numerical algorithms.
Hence, we will be concerned with the convergence and accuracy of our numerical methods not just in finding stationary states, but also in approximating the $L^2$ gradient descent dynamics in getting there, with an eye towards eventual application to grain boundary networks in the very general setting where each interface moves and may have a distinct anisotropic (normal dependent) surface tension.}

The surface tensions $\sigma_{LS}$, $\sigma_{VS}$, and $\sigma_{VL}$ of the interfaces $\Gamma_{LS}$, $\Gamma_{VS}$, and $\Gamma_{VL}$, respectively, will thus be allowed to be normal dependent (anisotropic).
In this two dimensional setting, we will write $\sigma(\theta)$ as shorthand for $\sigma(\cos\theta,\sin\theta)$.
The convexity assumption on the surface tensions, necessary for well-posedness of models such as (\ref{eq:multiphase}), can then be expressed as
\begin{equation}
\sigma''(\theta) + \sigma(\theta) \geq 0
\end{equation}
\textcolor{black}{which must hold for the surface tension associated with each interface in the network.
The condition is necessary for the lower semicontinuity of the model, without which minimizing sequences would develop fine scale oscillations; see e.g. \cite{morgan1997}.
It is also implied by our assumption (\ref{eq:convexsigma}).}
The boundary $\partial W_\sigma$ of the corresponding Wulff shape $W_\sigma$ can be described by the parameterized curve $(x(\theta),y(\theta))$, where
\begin{equation}
    \begin{split}
        &x(\theta) = -\sigma \left(\theta-\frac{\pi}{2}\right)\sin\theta - \sigma'\left(\theta-\frac{\pi}{2}\right)\cos\theta \; \mbox{, and}\\
        & y(\theta) = \sigma \left(\theta-\frac{\pi}{2}\right)\cos\theta - \sigma'\left(\theta-\frac{\pi}{2}\right)\sin\theta
    \end{split}
\end{equation}

\noindent The total interfacial energy we consider is the following special case of (\ref{eq:multiphase}):
\begin{equation}
\label{eq:sf_energy}
    \mathcal{E}(L) = \int_{\Gamma_{VL}} \sigma_{VL}(n(x)) ds(x) + \int_{\Gamma_{LS}} \sigma_{LS}(n(x)) ds(x) + \int_{\Gamma_{VS}} \sigma_{VS}((n(x)) ds(x) 
\end{equation}
\noindent Stationary state of this system is obtained by minimizing (\ref{eq:sf_energy}) under the constraint of \textcolor{black}{area} conservation:
\begin{equation}
    L_{eq} = \argmin_{|L| = A}  \mathcal{E}(L)
\end{equation}

\noindent For $x\in \Gamma_{VL}$, let $\theta(x)$ denote the angle that the unit outward normal vector of $\Gamma_{VL}$ makes with the positive $x$-axis.
As noted in \cite{bao_jiang_srolovitz_wang2017}, the stationary shape of the droplet satisfies: 
\begin{equation*}
\kappa(x)(\sigma_{VL}(\theta(x)) + \sigma''_{VL}(\theta(x))) = C
\end{equation*}
along the interface $\Gamma_{VL}$, and
\begin{equation}
\label{eq:contactcondition}
\sigma_{VL}(\theta(x^*))\cos\theta(x^*) -\sigma_{VL}'(\theta(x^*))\sin\theta(x^*) + \sigma_{LS}(x^*) - \sigma_{VS}(x^*) = 0
\end{equation}
at contact points (triple junctions) $x^* \in \Gamma_{VL} \cap \partial S$. 
Here, $\kappa(x)$ is the curvature of the curve $\Gamma_{VL}$ at the point $x$ and $C$ is a constant determined by the prescribed area $A$.


\noindent In this study, we are concerned with the $L^2$-gradient flow of energy (\ref{eq:sf_energy}) subject to the area constraint, which is known as area preserving weighted (anisotropic) curvature flow of an interface and a network.
The normal speed of the interface $\Gamma_{VL}$ under this flow is given by
\begin{equation}
\label{eq:normalspeed}
V_{\perp}(x) = m\big(\theta(x)\big)\bigg(\kappa\big(\theta(x)\big)\Big(\sigma_{VL}\big(\theta(x)\big) + \sigma''_{VL}\big(\theta(x)\big)\Big) + \mu\bigg)
\end{equation}
where $m:\textcolor{black}{\mathbb{S}^{1}}\rightarrow \mathbb{R}^+$ is a normal-dependent {\em mobility} factor, \textcolor{black}{and} $\mu$ is a Lagrangian multiplier that ensures area preservation.
Along with the junction condition (\ref{eq:contactcondition}), the normal speed (\ref{eq:normalspeed}) determines the evolution of the partition, at least for short time before topological changes or other singular behavior occurs.


\section{Threshold Dynamics for the VLS Model}
Threshold dynamics was proposed by Merriman, Bence, and Osher \cite{mbo92}, originally only for the isotropic equal surface tension case of model (\ref{eq:multiphase}), i.e. with $\sigma_{i,j}(n) \equiv 1$.
In the two phase setting, its extension to anisotropic surface tensions was considered in a number of earlier studies, e.g. \cite{ruuth_merriman, bbc, elsey_esedoglu_anisotropy}.
\textcolor{black}{In the multiphase setting, extension to unequal \cite{esedoglu_otto} and later to anisotropic  \cite{elsey_esedoglu_anisotropy,ejz} surface tensions relies on the variational formulation in \cite{esedoglu_otto} of the original equal, isotropic surface tension scheme from \cite{mbo92}.
(Also see \cite{laux_otto, laux_otto2} for rigorous results on convergence to the correct dynamics in the isotropic case).
These generalizatons of threshold dynamics generate a time discrete approximation to the evolution of the sets $\Sigma_j$ in the partition, starting from an initial condition $\Sigma_j^0$ as shown in Algorithm \ref{alg:matd}:}
\begin{algorithm}[H]
\caption{Multiphase Anisotropic Threshold Dynamics}
\label{alg:matd}
Using a time step size of $\delta t>0$, to obtain the sets $\Sigma_j^{k+1}$ at time $t=(k+1)\delta t$ from $\Sigma_j^k$ at time $t=k\delta t$:
\begin{algorithmic}[1]
\STATE \textbf{Convolution step}:
\begin{equation}
\label{eq:anisotropic_td1}
\psi_i^k =\sum_{j\not= i} K^{i,j}_{\delta t} * \mathbf{1}_{\Sigma_j^k}
\end{equation}
\STATE \textbf{Thresholding step}:
\begin{equation}
\label{eq:anisotropic_td2}
\Sigma_i^{k+1} = \{ x\in\Omega: \psi_i^k (x) = \min_{j} \psi_j^k(x) \}
\end{equation}
\end{algorithmic} 
\end{algorithm}
\noindent where $K^{i,j}$ are convolution kernels encoding the desired surface tensions $\sigma_{i,j}$ and mobilities $m_{i,j}$, and we write (for $\delta>0$):
\begin{equation}
K_\delta(x) = \frac{1}{\delta^d} K \left( \frac{x}{\delta} \right).
\end{equation}
Algorithm \ref{alg:matd} is derived by a systematic procedure from (and can be seen as implementing gradient descent for) the following nonlocal approximation to energy (\ref{eq:multiphase}):
\begin{equation}
\label{eq:nonlocalE}
E_{\delta t}(\Sigma_1,\Sigma_2,\ldots,\Sigma_N) = \textcolor{black}{\frac{1}{\sqrt{\delta t}}} \sum_{i\not =j} \int_{\Sigma_j} K^{i,j}_{\delta t} * \mathbf{1}_{\Sigma_i} \, dx
\end{equation}
\noindent  \textcolor{black}{Prior to this variational formulation, threshold dynamics schemes could not be correctly extended even to the 
isotropic, unequal surface tension case of model (\ref{eq:multiphase}) due to the formation of boundary layers of size $O(\sqrt{\delta t})$ at triple junctions.
The limiting interfacial dynamics the scheme converges to obeys an effective angle condition that differs from the apparent angle at the junction.
The variational formulation (\ref{eq:nonlocalE}) makes it possible to predict the effective angle in terms of the parameters of the scheme (e.g. choice of the convolution kernels and their weights), which would otherwise be very difficult.
The present study brings this powerful technique to level set methods, even in the presence of anisotropy, by exploting the precise connection between median filters and threshold dynamics noted in \cite{morel_book, esedoglu_guo_li}.
As a consequence of this subtlety, in verifying the convergence of the class of numerical methods studied here, one cannot measure the attainment of the angle conditions directly at a triple junction: the inner angle is simply irrelevant.
Instead, we assess the convergence of each phase (set) to the desired limiting motion e.g. in measure, which would necessarily imply the {\em effective angles} generated by the scheme are the correct Herring angles.}

In wetting / dewetting problems, we have three phases $V$ (vapor), $L$ (liquid) and $S$ (solid) partitioning the domain $\Omega$.
The solid phase $S$ is stationary, with a reasonably nice boundary $\partial S$, making the problem essentially two-phase: we can take $L$ to be the only variable, from which $V$ is determined essentially as $V=\Omega \setminus (L\cup S)$.
Nevertheless, the model features contact points (triple junctions) with interesting junction conditions, and is therefore a good test bed for multiphase curvature motion algorithms.
In addition, the \textcolor{black}{area (volume in higher dimensions)} of the liquid and vapor phases are preserved.
The energy in this context can be written as:
\begin{equation}
\label{eq:ELV}
    E_{\delta t}(L,V) = \textcolor{black}{\frac{1}{\sqrt{\delta t}}}\left( \int_L K^{VL}_{\delta t} *\mathbf{1}_V dx + \int_S K^{VS}_{\delta t} *\mathbf{1}_V dx +   \int_S K^{LS}_{\delta t} *\mathbf{1}_L dx \right)
\end{equation}
where $K^{VL}$, $K^{SL}$ and $K^{VS}$ are kernels that will be chosen to capture the anisotropies $\sigma_{VL}$, $\sigma_{LS}$, $\sigma_{VS}$.
Of course, since $| V \vartriangle (S^c \cap L) | = 0$, the variable $V$ can be eliminated from  (\ref{eq:ELV}), allowing the energy to be expressed even more simply as
\begin{equation}
\label{eq:EL}
E_{\delta t}(L) = \textcolor{black}{\frac{1}{\sqrt{\delta t}}} \int_L \left( K_{\delta t}^{LS} - K_{\delta t}^{VS} - K_{\delta t}^{VL} \right) * \mathbf{1}_S - K^{VL}_{\delta t} * \mathbf{1}_L \, dx
\end{equation}
up to a constant (taking into account that the volume of $L$ is preserved).
We will consider discrete in time, continuous in space schemes approximating the gradient descent dynamics associated with energy (\ref{eq:EL}) in the sense discussed in the previous section.
These schemes produce an approximation $L^{n+1}$ at time $t=(n+1)\delta t$ from the approximation $L^n$ at time $t=n\delta t$ by minimizing (\ref{eq:EL}) plus a movement limiter (proximal) term.
As in \cite{esedoglu_otto}, the resulting variational problems will be relaxed over functions satisfying a box constraint.
To that end, define the convex set $\mathcal{K}_1$ of admissible functions $u_L$ as follows:
\begin{equation}
\label{eq:K1}
\mathcal{K}_1 = \left\{ u_L \in BV(\Omega) \, : \, 0\leq u_L(x) \leq \mathbf{1}_{S^c}(x) \right\}
\end{equation}
and consider dissipating the following energy with respect to $u_L$ over $\mathcal{K}_1$
\begin{multline}
\label{eq:Eumu}
\tilde{E}_{\delta t} (u_L,\mu) = \textcolor{black}{\frac{1}{\sqrt{\delta t}}} {\textcolor{black}{\int_\Omega}} u_L \left( K_{\delta t}^{LS} - K_{\delta t}^{VS} - K_{\delta t}^{VL} \right) * \mathbf{1}_S\\
- u_L K^{VL}_{\delta t} * u_L + \mu \left( u_L - \frac{A}{|\Omega|} \right) \, dx
\end{multline}
for an appropriate choice of the Lagrange multiplier $\mu$.
Following \cite{esedoglu_otto}, a threshold dynamics algorithm can now be derived from (\ref{eq:Eumu}) by minimizing the linearization of (\ref{eq:Eumu}) at $u_L^n$
\begin{equation}
\mathcal{L}(v,u_L^n,\mu^n) =\textcolor{black}{\frac{1}{\sqrt{\delta t}}} {\textcolor{black}{\int_\Omega}} v \Bigg\{ \left( K_{\delta t}^{LS} - K_{\delta t}^{VS} - K_{\delta t}^{VL} \right) * \mathbf{1}_S
- 2 K^{VL}_{\delta t} * u_L^n+ \mu^n \Bigg\} \, dx
\end{equation}
with respect to $v \in \mathcal{K}_1$.
The result is easily seen to be the binary function
\begin{equation}
u_L^{n+1} = 
    \begin{cases}
    1 & \text{ if } K^{VL}_{\delta t} * (\mathbf{1}_{S^c}- 2 u^n_L) +(K^{LS}_{\delta t} - K^{VS}_{\delta t}) * \mathbf{1}_S  + \tilde{\mu}^n \leq 0\\
    0 & \text{ otherwise }
    \end{cases}
\end{equation}
on $S^c$, from which $L^{n+1}$ is determined; \textcolor{black}{here $\tilde{\mu}^n = \mu^n - \int K^{VL} \, dx$}. 
It turns out to be equivalent to minimizing movements for (\ref{eq:Eumu}) with a particular movement limiter term under the condition that $\widehat{K}^{VL} \geq 0$; see \cite{esedoglu_otto} for details.
The resulting scheme is described in Algorithm \ref{alg:pvs_td_anisotropic}.
\begin{algorithm}[H]
\caption{Threshold Dynamics for the VLS Model}
Find $\mu^n$ so that the following rule for obtaining $L^{n+1}$ from $L^n$ ensures $|L^{n+1}|=A$:
\label{alg:pvs_td_anisotropic}
\begin{algorithmic}[1]
\STATE Convolution step:
\begin{equation}
\label{eq:pvs_td1}
\psi^n = \Big(K^{VL}_{\delta t} * (\mathbf{1}_{S^c}- 2 \mathbf{1}_{L^n}) + (K^{LS}_{\delta t} - K^{VS}_{\delta t}) * \mathbf{1}_S  + \mu^n \Big)
\end{equation}
\STATE Thresholding step:
\begin{equation}
\label{eq:pvs_td2}
L^{n+1} = \{x \in S^c \, : \, \psi^n \leq 0 \}
\end{equation}
\end{algorithmic}
\end{algorithm}
\noindent The one-dimensional search for the scalar Lagrange multiplier $\mu^n$ at every time step can be accomplished in many ways, e.g. via bisection, the secant method, or a combination of the two.
Let us introduce notation for applying one step of this scheme:
\begin{equation}
L^{n+1} = \mathbf{T}_{\delta t} (L^n,V^n,S,\mu^n).
\end{equation}

\section{Median Filter for the VLS Model}
\label{sec:medianfilter}
The novelty of the present study is a {\em level set method} in the form of a {\em vector valued median filter} that will be derived from the threshold dynamics Algorithm \ref{alg:pvs_td_anisotropic} following \cite{esedoglu_guo_li} (which in turn builds up on \cite{oberman2004}).
A potential advantage of the level set method over threshold dynamics is its representation of interfaces using regular (vs. characteristic) functions, allowing for subgrid accuracy in locating the interface via interpolation on uniform spatial grids; lack of any subgrid accuracy in naive implementations of threshold dynamics can result in interfaces getting pinned (but see \cite{ruuththesis, ruuth0, ruuth1} for efficient implementation of threshold dynamics on adaptive grids as another solution to this drawback).

With that in mind, let $\phi_L^n(x)$, $\phi_V^n(x)$, $\phi_S(x)$ denote the level set functions the zero super-level sets of which represent the liquid, vapor, and solid phases.
Following \cite{esedoglu_guo_li}, our goal is to \textcolor{black}{give} an algorithm for updating $\phi^n_L$ and $\phi^n_V$ to $\phi^{n+1}_L$ and $\phi^{n+1}_V$ that is consistent with Algorithm \ref{alg:pvs_td_anisotropic}, i.e.
\begin{equation}
L^n = \{ x \, : \, \phi_L^n(x) \geq 0 \} \; , \;  V^n = \{ x \, : \, \phi_V^n(x) \geq 0 \} \mbox{ and } S = \{ x \, : \, \phi_S(x) \geq 0 \}.
\end{equation}
for all $n$, provided that $\{ \phi^0_L \geq 0 \} = L^0$, $\{ \phi^0_V \geq 0\} = V^0$, and $\{ \phi_S \geq 0 \} = S$ where $V^n$, $L^n$, and $S$ as before denote the liquid and vapor phases at time step $n$ generated by Algorithm \ref{alg:pvs_td_anisotropic}.
In other words, the zero-level sets of $\phi^n_L$, and $\phi^n_V$ should implement exactly Algorithm \ref{alg:pvs_td_anisotropic}.
At the same time, Algorithm \ref{alg:pvs_td_anisotropic} should be extended to other level sets of $\phi^n_V$ and $\phi^n_L$, in such a way that $\phi^n_L$, and $\phi^n_V$ and remain regular if $\phi^0_L$, and $\phi^0_V$ are.

To that end, as in \cite{esedoglu_guo_li}, define the partition at level $\lambda$ and time step $n$ {\em from the point of view of the liquid phase} as:
\begin{equation}
\label{eq:povL}
\left\{
\begin{split}
L^n(\lambda) &= \{ x \, : \, \phi^n_L(x) \geq \lambda \},\\
V_L^n(\lambda) &= \{ x \, : \, \phi^n_L(x) < \lambda \mbox{ and } \phi^n_V(x) \geq \phi_S(x) \},\\
S_L(\lambda) &= \{ x \, : \, \phi^n_L(x) < \lambda \mbox{ and } \phi_S(x) > \phi^n_V(x) \},\\
\end{split}
\right.
\end{equation}
and the partition at level $\lambda$ and time step $n$ {\em from the point of view of the vapor phase} as:
\begin{equation}
\label{eq:povV}
\left\{
\begin{split}
V^n(\lambda) &= \{ x \, : \, \phi^n_V(x) \geq \lambda \},\\
L_V^n(\lambda) &= \{ x \, : \, \phi^n_V(x) < \lambda \mbox{ and } \phi^n_L(x) \geq \phi_S(x) \},\\
S_V(\lambda) &= \{ x \, : \, \phi^n_V(x) < \lambda \mbox{ and } \phi_S(x) > \phi^n_L(x) \},\\
\end{split}
\right.
\end{equation}
The $\lambda$-super level sets of $\phi^n_L$ will be updated by applying thresholding Algorithm \ref{alg:pvs_td_anisotropic} to the partition (\ref{eq:povL}), whereas the $\lambda$-super level sets of $\phi^n_V$ will be updated by applying the same thresholding algorithm to the partition (\ref{eq:povV}), using the same ($\lambda$ independent) Lagrange multiplier, ensuring volume preservation only for the $0$-super level sets:
\begin{equation}
\textcolor{black}{
\label{eq:levels}
\begin{split}
L^{n+1}(\lambda) &= \mathbf{T}_{\delta t} (L^n(\lambda),V_L^n(\lambda),S_L(\lambda),\mu^n)\\
= \Big\{ x& \in S_L^c(\lambda) \, : \, K^{VL}_{\delta t} * (\mathbf{1}_{S_L^c(\lambda)}- 2 \mathbf{1}_{L^n(\lambda)}) + (K^{LS}_{\delta t} - K^{VS}_{\delta t}) * \mathbf{1}_{S_L(\lambda)}  + \mu^n \leq 0 \Big\}\\
V^{n+1}(\lambda) &= \mathbf{T}_{\delta t} (V^n(\lambda),L_V^n(\lambda),S_V(\lambda),\mu^n)\\
= \Big\{ x& \in S_V^c(\lambda) \, : \, K^{VL}_{\delta t} * (\mathbf{1}_{S_V^c(\lambda)}- 2 \mathbf{1}_{V^n(\lambda)}) + (K^{VS}_{\delta t} - K^{LS}_{\delta t}) * \mathbf{1}_{S_V(\lambda)}  + \mu^n \leq 0 \Big\}\\
\end{split}}
\end{equation}
We note the following simple comparison property which ensures there are functions $\phi^{n+1}_L$ and $\phi^{n+1}_V$ having $L^{n+1}(\lambda)$ and $V^{n+1}(\lambda)$ from equation (\ref{eq:levels}) as their $\lambda$-super level sets, respectively, though this requires a restrictive condition on the kernels:
\begin{claim}
\label{claim:max}
Assume the pointwise triangle inequality $K^{LS} + K^{VL} - K^{VS}\geq 0$ and $K^{VS} + K^{VL} - K^{LS}\geq 0$ on the kernels.
If $\lambda_1 > \lambda_2$, then
$$ L^{n+1}(\lambda_1) \subseteq L^{n+1}(\lambda_2) \mbox{ and } V^{n+1}(\lambda_1) \subseteq V^{n+1}(\lambda_2) $$
where $L^{n+1}$ and $V^{n+1}$ are given by (\ref{eq:levels}).
\end{claim}

\begin{proof}
\textcolor{black}{
We will show $L^{n+1}(\lambda_1) \subseteq L^{n+1}(\lambda_2) $;
the proof of $V^{n+1}(\lambda_1) \subseteq V^{n+1}(\lambda_2)$ is identical.
To that end, let $\gamma_\lambda = \left( K^{LS} - K^{VS} - K^{VL} \right) * \mathbf{1}_{S_L(\lambda)}
- 2 K^{VL} * \mathbf{1}_{L(\lambda)}$ and $\Tilde{K} =  K^{LS} - K^{VS} - K^{VL}$. 
We check
\begin{center}
 $\lambda_1 \geq \lambda_2$ implies $\gamma_{\lambda_1} \geq \gamma_{\lambda_2}$. 
\end{center}
Indeed,
\begin{equation}
\label{eq:maxproof}
\begin{split}
       \gamma_{\lambda_1} - \gamma_{\lambda_2} &= \Tilde{K} * \left(\mathbf{1}_{S_L(\lambda_1)} - \mathbf{1}_{S_L(\lambda_2)}\right) - 2 K^{VL} * \left(\mathbf{1}_{L(\lambda_1)} - \mathbf{1}_{L(\lambda_2)}\right)\\
      & = \Tilde{K} * \mathbf{1}_{S_L(\lambda_1)\setminus S_L(\lambda_2)} + 2 K^{VL} * \mathbf{1}_{L(\lambda_2)\setminus L(\lambda_1)}\\
      & \geq \Tilde{K} * \mathbf{1}_{S_L(\lambda_1)\setminus S_L(\lambda_2)} + 2 K^{VL} * \mathbf{1}_{S_L(\lambda_1)\setminus S_L(\lambda_2)}\\
      & =  \left(K^{LS} + K^{VL} - K^{VS}\right) * \mathbf{1}_{S_L(\lambda_1)\setminus S_L(\lambda_2)}\\
      & \geq 0
\end{split}
\end{equation}
The second equality holds because $S_L(\lambda_2) \subseteq S_L(\lambda_1)$, and $L(\lambda_1) \subseteq \textcolor{black}{L}(\lambda_2)$, the third inequlity holds because $S_L(\lambda_2) \setminus S_L(\lambda_1) \subseteq L(\lambda_1) \setminus L(\lambda_2)$, and the last inequality holds because $K^{LS} + K^{VL} - K^{VS} \geq 0$, verifying $L^{n+1}(\lambda_1) \subseteq L^{n+1}(\lambda_2)$. 
}
\end{proof}

\begin{remark}
\textcolor{black}{According to formulas (17) and (41) of \cite{ejz}},
the pointwise triangle inequality required of the kernels in the hypothesis of Claim \ref{claim:max} implies the following particularly strong assumption on the surface tensions:
\begin{equation}
\label{eq:st}
 (\sigma+\sigma'')_{i,j} + (\sigma+\sigma'')_{j,k} \geq (\sigma+\sigma'')_{i,k} 
\end{equation}
whenever $i$, $j$, and $k$ are distinct.
Several examples in Section \ref{sec:numerical_prescirbed_dynamics} violate this undesirably strong assumption, but we see no adverse effects in the numerical results.
\end{remark}

The crucial observation \textcolor{black}{(see Lemma 3.1 in 
 \cite{esedoglu_guo_li})}, is that all this (infinitely many applications of threshold dynamics algorithm) can be accomplished concurrently (on all super level sets) by a vectorial local median filter acting on the vector valued function $\phi^n_V(x), (\phi_L^n(x), \phi_S^n(x))$.
The resulting numerical scheme is described in Algorithm \ref{alg:pvs_mf_sorting1}.
Its most significant step is the essential ingredient in any local median filter: {\em sorting} of (in this case the level set) function values. \textcolor{black}{Choose $m$ points $y_1,...,y_m$ that sample approximately uniformly the circle $\partial B_r(0)$, with $r = \sqrt{\delta t}$:} 

\begin{algorithm}[H]
\caption{Median Filter for the VLS Model}
\label{alg:pvs_mf_sorting1}
Find $\mu^n$ so that the following update for getting $\phi^{n+1}_L$ from $\phi^n_L$ ensures $|\{ \phi^{n+1}_L \geq 0 \}| = A$:
\begin{algorithmic}[1]
\STATE \textbf{Sort}: Choose $m$ neighborhood points around $x$ and sort the level set values $\{\phi_L^n(x+y_1), \phi_L^n(x+y_2),...,\phi_L^n(x+y_m)\}$ so that the permutation $p: \{1,2,...,m\} \rightarrow \{1,2,...,m\}$ satisfies
\begin{equation}
\label{eq:mf1}
\phi_L^n(x+y_{p(1)}) \leq \phi_L^n(x+y_{p(2)}) \leq \cdots \leq \phi_L^n(x+y_{p(m)})
\end{equation}
\STATE \textbf{Set}: $C_1 = 0$, $C_2 = ||K^{VL}||$ and $\ell = 1$.
\STATE \textbf{While} $C_1\leq C_2 + \mu^n$ do 
\STATE \text{~~~~} \textbf{If} $\phi^n_V(x+y_{p(\ell)}) \geq \phi^n_S(x+y_{p(\ell)})$
\STATE \text{~~~~~~~~} $C_1\leftarrow C_1+ K^{VL}(y_{p(\ell)})$ and $C_2\leftarrow C_2 - K^{VL}(y_{p(\ell)})$
\STATE \text{~~~~} \textbf{Else}
\STATE \text{~~~~~~~~} $C_1\leftarrow C_1+ K^{LS}(y_{p(\ell)})$ and $C_2\leftarrow C_2 + K^{VS}(y_{p(\ell)}) - K^{VL}(y_{p(\ell)})$
\STATE  \text{~~~~} $\ell \leftarrow \ell+1$
\STATE \textbf{End While}
\STATE $\phi_L^{n+1}\textcolor{black}{(x)} = \frac{1}{2} \Big( \phi_L^n(x+y_{p(\ell-1)}) + \phi_L^n(x+y_{p(\ell)})\Big)$.
\STATE \textbf{Set} $\phi_L^{n+1}\textcolor{black}{(x)} \leftarrow \min\{\phi_L^{n+1}\textcolor{black}{(x)}, -\phi_S\textcolor{black}{(x)}\}$.
\end{algorithmic}
Update $\phi^n_V$ to $\phi^{n+1}_V$ analogously, \textcolor{black}{by interchanging the $L$ and $V$ in superscripts and superscripts from the above algorithm.}
\end{algorithm}

\begin{center}
    \begin{figure}
        \centering
        \includegraphics[scale = 0.75]{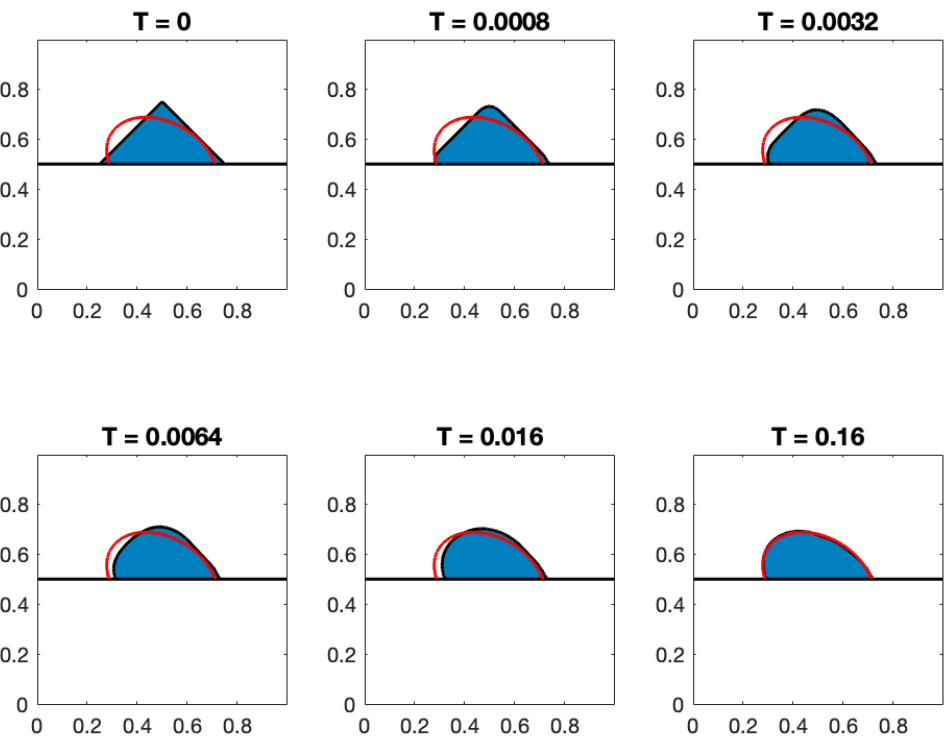}
        \caption{Dynamics of the anisotropic droplet using the sorting algorithm \ref{alg:pvs_mf_sorting1}. The filled component is the liquid droplet at different times $T$, the red line is the theoretic Wulff shape associated with the surface tension, and the black lines are interfaces.}
        \label{fig:sorting_dynamics}
    \end{figure}
\end{center}

\textcolor{black}{At every time step the value of the Lagrange multiplier $\mu^n$ is determined via bisection. By construction, under Algorithm \ref{alg:pvs_mf_sorting1}, the zero super level sets $\{\phi^n_L\geq0\}$ , $\{\phi^n_V\geq0\}$, and $\{\phi^n_S\geq 0\}$ of the three phases are updated by a discrete implementation of Algorithm \ref{alg:pvs_td_anisotropic}, and hence form a partition.}
Naive implementation of median filters by sampling and sorting $m$ points from every neighborhood typically requires a large $m$ for accurate results.
Figure \ref{fig:sorting_dynamics} illustrates a numerical result using the above sorting-based algorithm for a liquid droplet with the surface tensions \textcolor{black}{$\sigma_{VL}(\theta) = \sqrt{1+\sin^2(\theta-\frac{\pi}{3})}$}, $\sigma_{SL} = \sigma_{SV} = 1$.
The computation was performed using $400\times400$ grid points, \textcolor{black}{and $100$} points sampled from every circular neighborhood, and a time step size of $\delta t=0.0004$.
\textcolor{black}{We use the $L^1$ norm  (an approximation to the area of symmetric difference between the shapes) and the $L^\infty$ norm (an approximation to the Hausdorff distance between their boundaries) to measure the error.
The $L^1$ error in the computed Wulff shape at final time $T = 0.16$ is $0.0031$; the $L^\infty$ error is $0.0085$}.
\textcolor{black}{An efficient implementation}, using an extension of the accurate median filter described in \cite{esedoglu_guo_li} to general anisotropic surface tensions and mobilities, was used in generating the numerical convergence studies presented in subsequent sections.

\section{Anisotropic Kernels for Median Filters}
\label{sec:Conv_kernel_construction}

\noindent A number of earlier studies discuss how to find a convolution kernel $K$ for a desired anisotropic (normal dependent) surface tension and anisotropic mobility pair $(\sigma_*,\textcolor{black}{m}_*)$ to be used in threshold dynamics algorithms; see e.g. \cite{bbc, elsey_esedoglu_anisotropy, esedoglu_jacobs, ejz}.
In particular, \cite{ejz} gives an explicit construction that in two dimensions guarantees the positivity of $K$ for essentially all anisotropic surface tensions (provided the underlying surface energy model is well-posed) and a broad range of anisotropic mobilities.
(Note, however, \cite{elsey_esedoglu_anisotropy, ejz} also reveal obstructions to positivity in dimensions three and higher).

For the extension of these algorithms to level set methods in the form of median filters as discussed in \cite{esedoglu_guo_li} and Section \ref{sec:medianfilter} of the present study, positivity of the kernels are essential.
Moreover, the highly accurate implementation of median filters discussed in \cite{esedoglu_guo_li} favors kernels concentrated on a few circles (diffuse kernels can be accommodated, but require another layer of approximation).
Therefore, in this section, we revisit the positive kernel construction of \cite{ejz}, but in the special context of kernels concentrated on a few circles.

Given a desired, possibly anisotropic pair $(\sigma_*,\textcolor{black}{m}_*)$ of surface tension and mobility, a corresponding kernel $K$ can be found by solving the following system of integral equations from \cite{ejz}:
\begin{equation}
\label{eq:inteq1}
    \begin{cases}
        \int_{\mathbb{R}^d} K(x)|n\cdot x|dx = \sigma_*(n)\\
        \int_{n^{\perp}}K(x)dH^{d-1}(x) = \frac{1}{\textcolor{black}{m}_*(n)}
    \end{cases}
\end{equation}
System (\ref{eq:inteq1}) is more conveniently expressed using the cosine and spherical Radon transforms, $\mathcal{T}$ and $\mathcal{J}_s$, respectively:
\begin{equation}
\label{eq:inteq2}
    \begin{cases}
        \int_{0}^{\infty} K(rn)r^d dr = \mathcal{T}^{-1}\sigma_*(n)\\
        \int_{0}^{\infty}K(rn) r^{d-2}dr = \mathcal{J}_s^{-1}[\frac{1}{\textcolor{black}{m}_*}](n)
    \end{cases}
\end{equation}
The solution, when it exists, is most definitely not unique; hence the many different kernel constructions in existing literature.
This abundance makes it possible to impose further conditions, such as positivity, on the kernel in dimension $d=2$.
There, system (\ref{eq:inteq2}) can be further simplified by noting
\begin{equation}
\label{eq:sf_construt}
    \mathcal{T}^{-1}\sigma_*(\theta) = \frac{1}{4}\Big(\sigma_*''(\theta-\frac{\pi}{2}) + \sigma_*(\theta-\frac{\pi}{2})\Big)
\end{equation}
and
\begin{equation}
\label{eq:mobility_construt}
    \mathcal{J}_s^{-1}\Big[\frac{1}{\textcolor{black}{m}_*}\Big](\theta) = \frac{1}{\textcolor{black}{m}_*(\theta-\frac{\pi}{2})}
\end{equation}
Here, we claim that we can further require the kernel to be supported on two concentric circles of differing radii:

\begin{claim}
For any pair of $(\sigma_*, \textcolor{black}{m}_*)$ with $\sigma''_*  + \sigma_* > 0$ and $\textcolor{black}{m}_*>0$, there exists a {\em positive} kernel $K$ {\em supported on two concentric circles} having $\sigma_*$ and $\textcolor{black}{m}_*$ as its corresponding surface tension and mobility, as given in (\ref{eq:inteq1}).
\end{claim}  

\begin{proof}
    Write $K(r,\theta) = \omega_1(\theta)\delta(r-R_1) + \omega_2(\theta)\delta(r-R_2)$, where $\omega_i(\theta) > 0$ for any $\theta$, $i = 1,2$. Then $\omega_i$'s and $R_i$'s satisfy the following equation:
    \begin{equation}
        \begin{pmatrix}
            R_1^2 & R_2^2 \\
            1 & 1
        \end{pmatrix}
     \begin{pmatrix}
           \omega_1(\theta) \\
           \omega_2(\theta)
        \end{pmatrix} = 
        \begin{pmatrix}
           \frac{1}{4}\Big(\sigma_*(\theta-\frac{\pi}{2}) + \sigma_*''(\theta-\frac{\pi}{2})\Big) \\
         \frac{1}{\textcolor{black}{m}_*(\theta-\frac{\pi}{2})}
        \end{pmatrix}
    \end{equation}
    If the matrix on the \textcolor{black}{left}-hand side is invertible, we can then write
    \begin{equation}
   \label{eq:weight_ftn}
        \begin{pmatrix}
           \omega_1(\theta) \\
           \omega_2(\theta)
        \end{pmatrix} = \frac{1}{R_1^2 - R_2^2} \begin{pmatrix}
            1 & -R_2^2 \\
            -1 & R_1^2
        \end{pmatrix}
        \begin{pmatrix}
           \frac{1}{4}\Big(\sigma_*(\theta-\frac{\pi}{2}) + \sigma_*''(\theta-\frac{\pi}{2})\Big) \\
         \frac{1}{\textcolor{black}{m}_*(\theta-\frac{\pi}{2})}
        \end{pmatrix}
    \end{equation}
    Given that $\sigma_* + \sigma_*'' > 0$ and $\textcolor{black}{m}_* >0$, as long as $R_2 \ll 1 \ll R_1$, the positivity of $\omega_i$'s are guranteed. More precisely, as long as $\max\{R_1,R_2\} > \max\{\frac{1}{4}\mu_*(\theta-\frac{\pi}{2})\Big(\sigma_*(\theta-\frac{\pi}{2}) + \sigma_*''(\theta-\frac{\pi}{2})\Big) \}^{\textcolor{black}{\frac{1}{2}}}$, $\min\{R_1,R_2\} < \min\{\frac{1}{4}\textcolor{black}{m}_*(\theta-\frac{\pi}{2})\Big(\sigma_*(\theta-\frac{\pi}{2}) + \sigma_*''(\theta-\frac{\pi}{2})\Big) \}^{\textcolor{black}{\frac{1}{2}}}$, the weights on the circles are always positive.
\end{proof}

This is the construction that will be used in the remaining of the paper, in particular in the numerical experiments and convergence studies in subsequent sections.
\textcolor{black}{For example, for the $4$-fold anisotropy $\sigma(\theta) = 1 - \frac{1}{16}\cos 4\theta$ and constant mobility $m(\theta) = 1$, formula (\ref{eq:weight_ftn}) gives
\begin{equation}
\begin{split}
    &\omega_1(\theta) = \frac{16-64 R_2^2+15\cos(4\theta)}{64(R_1^2-R_2^2)}\\
    &\omega_2(\theta) = \frac{-\frac{1}{4}+R_1^2 - \frac{15}{64}\cos(4\theta)}{R_1^2-R_2^2}\\
\end{split}
\end{equation}}

\section{Numerical Experiments}
\label{sec:numerical_experiment}

\noindent This section presents a variety of numerical results of Algorithm \ref{alg:pvs_mf_sorting1}  using the fast and accurate median computation techniques as proposed in \cite{esedoglu_guo_li}.
\textcolor{black}{We recall from \cite{esedoglu_guo_li} that in two dimensions, the multiphase median filter schemes such as \ref{alg:pvs_mf_sorting1} can be implemented at a computational cost comparable to the threshold dynamics version, but with the enhanced spatial accuracy afforded by the level set representation.
Just like threshold dynamics, multiphase median filters extend seamlessly to higher dimensions, but their efficient implementation in that setting remains a goal for the future.}

The accuracy in time of the dynamics is assessed by comparing against front tracking simulations with extremely fine finite differences discretization (forward Euler in time, and centered differences in space). In Section \ref{sec:numerical_anisotropic}, we test the convergence rate (in time) of a droplet with an anisotropic surface tension $\sigma_{VL}(\mathbf{n})$ but constant mobility, and constant surface tensions $\sigma_{VS}$ and $\sigma_{LS}$ without topological changes.
The initial shapes of the liquid droplets are triangles. In this example, the solid-nonsolid interface is \textcolor{black}{flat (i.e. the outward unit normal vector is fixed)}, and the convolution kernel is supported on a single circle with its weight function \textcolor{black}{$\omega$} determined by the surface tension.
In Section \ref{sec:numerical_prescirbed_dynamics}, we use the approach introduced in Section \ref{sec:Conv_kernel_construction} to construct kernels using prescribed surface tensions and non-trivial mobilities.
We have studied the accuracy (in time) of the dynamics which converges to the same Wulff shape but with different mobilities.
In Section \ref{sec:numerical_fully_anisotropic}, we extend Algorithm \ref{alg:pvs_mf_sorting1} to a fully anisotropic setting, i.e., $\sigma_{VS}$, $\sigma_{VL}$, and $\sigma_{LS}$ are all normal-dependent, the mobilities are non-trivial, and then we study the convergence rate (in time) of the Algorithm \ref{alg:pvs_mf_sorting1}.
In Section \ref{sec:numerical_topological_changes}, we present two examples of simulations involving topological changes, specifically of (wetting and dewetting) droplets moving on a substrate with a sinusoidal \textcolor{black}{substrate} interface and anisotropic surface tensions at the interface $\Gamma_{VL}$.

In all the examples, the computational domain is a discretization of $[0, 1]\times[0, 1]$.
\textcolor{black}{The convergence rates observed in all cases are on the whole consistent with the expected rate of $\frac{1}{2}$, which is the rate observed in threshold dynamics simulations in the presence of junctions, see e.g. \cite{ruuththesis}; this is because median filters are consistent with threshold dynamics on zero super level sets.
In some of our experiments, a slow down in the rate of convergence can be seen at the finest resolutions.
We expect this is a consequence of the very rapid decrease of error in the first few refinements.}

\subsection{Prescribed anisotropic $\sigma_{VL}$}
\label{sec:numerical_anisotropic}

\begin{figure}
  \begin{minipage}[b]{.4\linewidth}
    \centering
     \includegraphics[scale = 0.3]{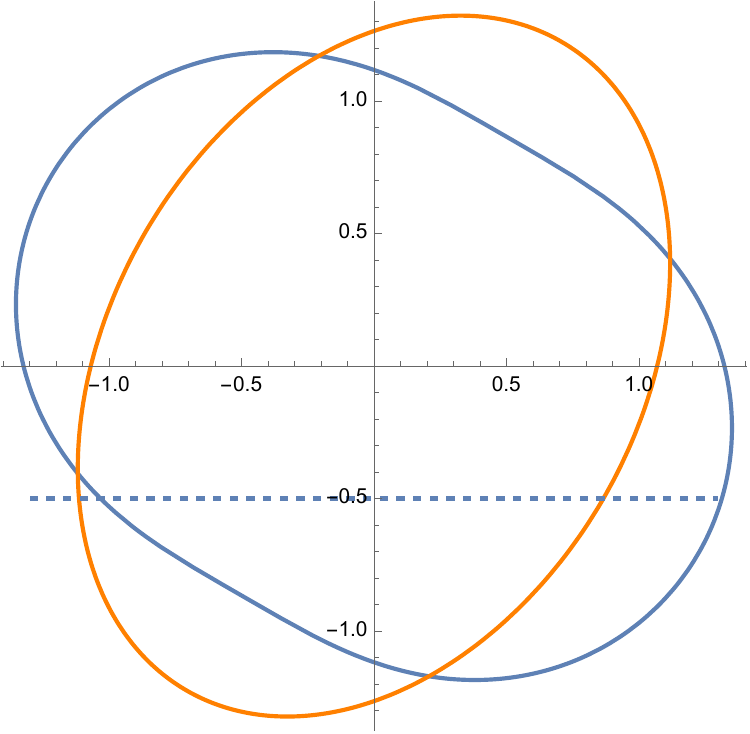}
\label{fig:anisotropic_one_circle_T0p008_wulff}
  \end{minipage}\hfill
  \begin{minipage}[b]{.5\linewidth}
  \centering
    \includegraphics[scale = 0.45]{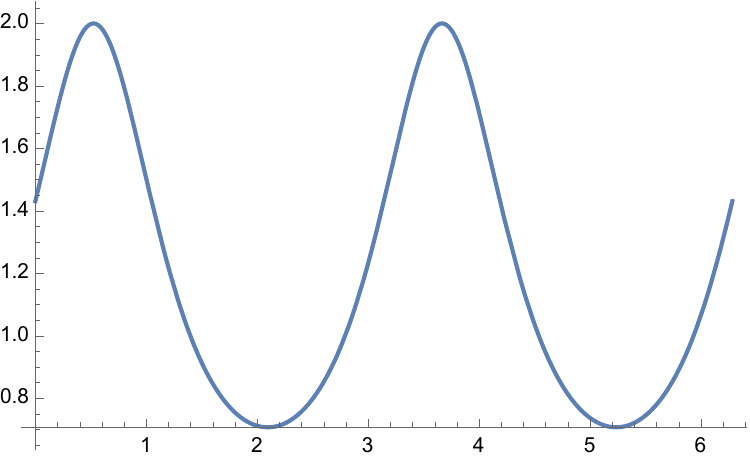}
\label{fig:anisotropic_one_circle_T0p008_weight_ftns}
  \end{minipage}
  \caption{Left: Polar plot for the surface tension $\sigma_{VL}$ and its Wulff shape. The blue line represents the plot for surface tension, and the orange line depicts the Wulff shape. The region of the Wulff shape above the dashed line (where $\sigma_{VS} - \sigma_{LS} = -0.5$) corresponds to the equilibrium state of this droplet.
Right: The weight function associated with the surface tension $\sigma_{VL}$ and $\textcolor{black}{m_{VL}}$.}
\label{fig:anisotropic_onecircle_info}
\end{figure}

In this section, we apply Algorithm \ref{alg:pvs_mf_sorting1} to an anisotropic example using a kernel supported on a single circle. The example has a flat solid-liquid interface. We chose the initial shape of the droplet to be an isosceles triangle with an area of 0.0625. In this experiment, we use a kernel supported on a single circle. The kernel was chosen to correspond to the desired surface tension $\sigma_{VL}(\theta) = \sqrt{1+\cos^2(\theta-\frac{\pi}{3})}$; however, the mobility is then determined to be $\textcolor{black}{m}_{VL}(\theta) = \frac{1}{\sigma_{VL}(\theta) + \sigma_{VL}''(\theta)} = \frac{(3 - \cos(\frac{\pi}{6} + 2\theta))^{3/2}}{4\sqrt{2}}$. The other two surface tensions were constants: $\sigma_{LS} = 1.5$ and $\sigma_{VS} = 1$. The Wulff shape associated with the surface tensions and the weight function determined by the surface tension and mobility is given in Figure \ref{fig:anisotropic_onecircle_info}. Initially, the triangle droplet was located on the solid surface at $y = 0.5$.\\

\begin{figure}
    \centering
    \includegraphics[width=0.45\linewidth]{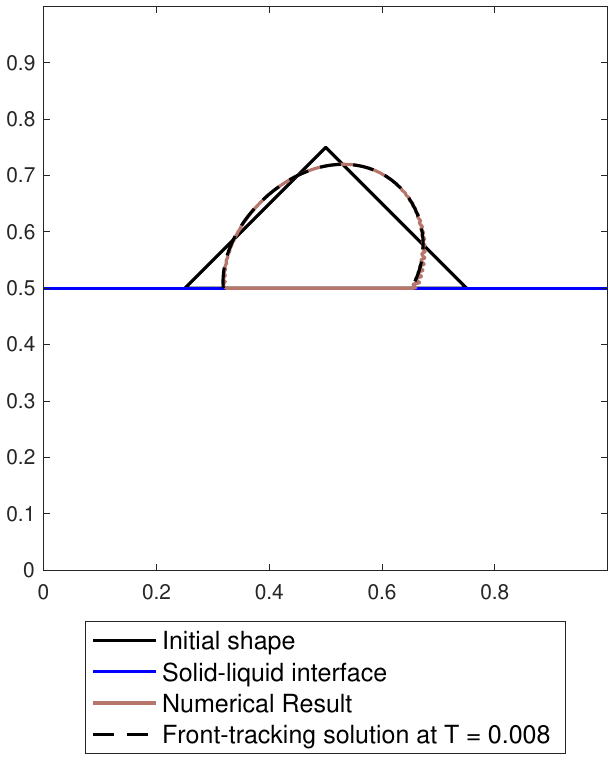}
    \caption{Anisotropic dewetting experiment with a droplet on a solid substrate. The black solid line is the initial liquid-vapor interface, the blue solid line is the solid surface, the black dashed line is the front-tracking numerical results with extremely fine $dx$ and $dt$, the pink solid line is the numerical result by median filter scheme using the finest $dt$ and $dx$ in Table \ref{table:one_circle_convergence}.}
    \label{fig:anisotropic_one_and_two_examples}
\end{figure}

\begin{figure}[H]
  \begin{minipage}[b]{.45\linewidth}
    \centering
    \includegraphics[scale=0.45]{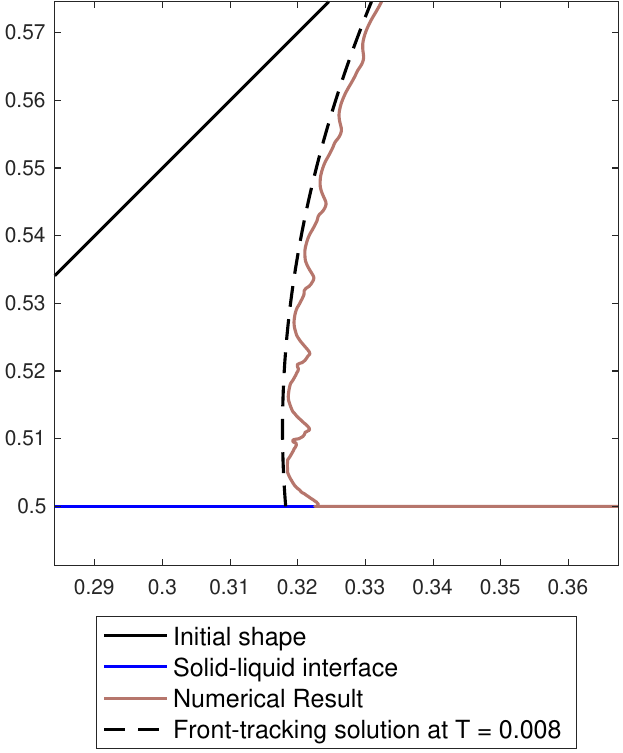}
  \end{minipage}\hfill
  \begin{minipage}[b]{.45\linewidth}
    \centering
    \includegraphics[scale=0.45]{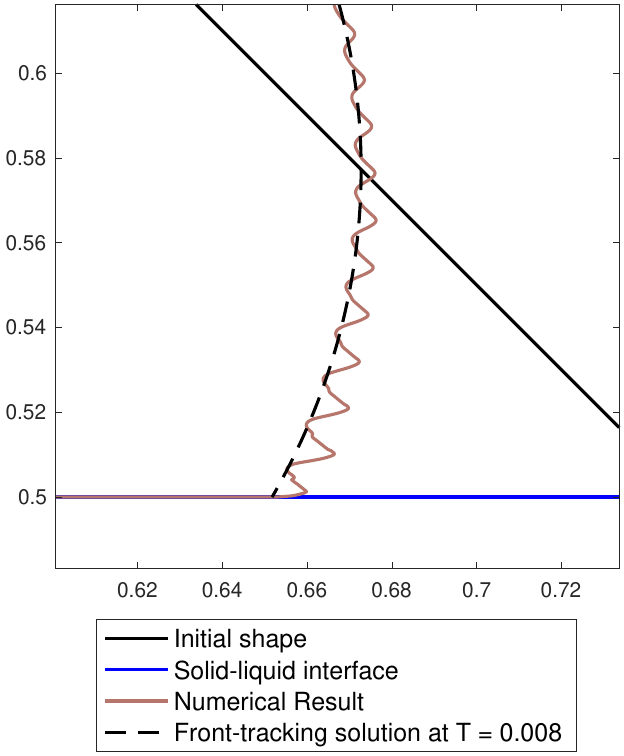}
  \end{minipage}
  \caption{Left and right junctions. Oscillations are observable, but do not appear to inhibit convergence; see error Table \ref{table:one_circle_convergence}.}
\label{fig:anisotropic_one_circle_T0p008_angles}
\end{figure}

\begin{figure}
  \begin{minipage}[b]{.45\linewidth}
    \centering
    \begin{adjustbox}{width=\columnwidth,center}  
    \begin{tabular}[b]{c|c|c|c} 
    \hline
     Number of time steps & $1/dx$ & $L^{1}$ Error & Order \\
     \hline

     \rule{0pt}{15pt}2 & 50 &  0.031213    & \textcolor{black}{-}\\
     \hline
     \rule{0pt}{15pt}4 & 100 &  0.007801     & 2.00  \\
     \hline
    \rule{0pt}{15pt} 8 & 200 & 0.0043854      & 0.83 \\
     \hline
     \rule{0pt}{15pt}16 & 400 & 0.0030301   & 0.53  \\
     \hline
     \rule{0pt}{15pt}32 & 800 & 0.0019095    & 0.67 \\
     \hline
     \rule{0pt}{15pt}64 & 1600 & 0.0011307   & 0.76  \\
     \hline
    \rule{0pt}{15pt} 128 & 3200 &  0.0008001  &  0.50   
\end{tabular}
    \end{adjustbox}
    \captionof{table}{Error table for the dynamic convergence test of a droplet with anisotropic surface tension $\sigma_{VL}(\theta) = \sqrt{1+\cos^2(\theta-\frac{\pi}{3})}$ and mobility $\textcolor{black}{m}_{VL}(\theta) = \frac{(3 - \cos(\frac{\pi}{3} + 2\theta))^{3/2}}{4\sqrt{2}}$ on a curved solid substrate. Final time $T = 0.008$.}
    \label{table:one_circle_convergence}
  \end{minipage}\hfill
  \begin{minipage}[b]{.45\linewidth}
  \centering
    \includegraphics[scale = 0.32]{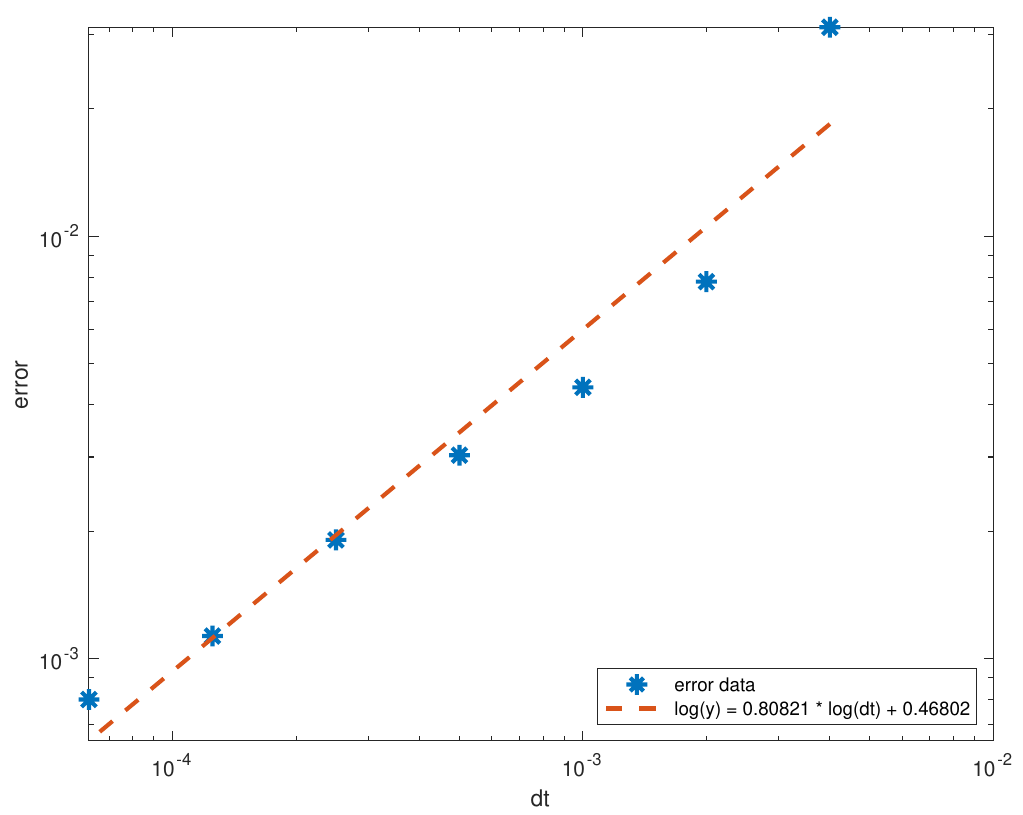}
\label{fig:anisotropic_one_circle_T0p008_error}
\captionof{figure}{The error data with fitting line. The slope of the fitting line in the log-log plot is around 0.80821. The expected convergence rate is $\frac{1}{2}$.}
\label{fig:anisotropic_onecircle_error}
  \end{minipage}
\end{figure}

\noindent In Table \ref{table:one_circle_convergence} and \textcolor{black}{Figure \ref{fig:anisotropic_onecircle_error}}, we present the corresponding $L^1$ errors and convergence rates in time. 
The final time $T$ is set to be 0.008. As noted in \cite{esedoglu_guo_li}, the condition $\hat{K} \geq 0$ is not satisfied by the weight function supported on a circle. This suggests that the energy stability of the scheme may be violated. It is a possible explanation for the slight oscillations observed around the junctions in Figures \ref{fig:anisotropic_one_circle_T0p008_angles}.\\

\noindent The significant negative part of the Fourier transform of the kernel concentrated on one circle, particularly around the origin, is possibly the reason for the more pronounced oscillations observed in the experiment, \textcolor{black}{as shown in Figures \ref{fig:anisotropic_one_and_two_examples} and \ref{fig:anisotropic_one_circle_T0p008_angles}.} However, \textcolor{black}{it might be possible to resolve this problem} by adapting the slightly "more Gauss-Seidel" version of the threshold dynamics algorithm, as described by \cite{esedoglu_jacobs}, for wetting and dewetting phenomena to its corresponding median filter scheme, which requires only $K \geq 0$ for energy stability, as pointed out in \cite{esedoglu_guo_li}.


\subsection{Prescribed anisotropic $\sigma_{VL}$ and $\textcolor{black}{m}_{VL}$}
\label{sec:numerical_prescirbed_dynamics}
\begin{figure}
 \begin{minipage}[b]{.4\linewidth}
    \centering
    \includegraphics[scale= 0.3]{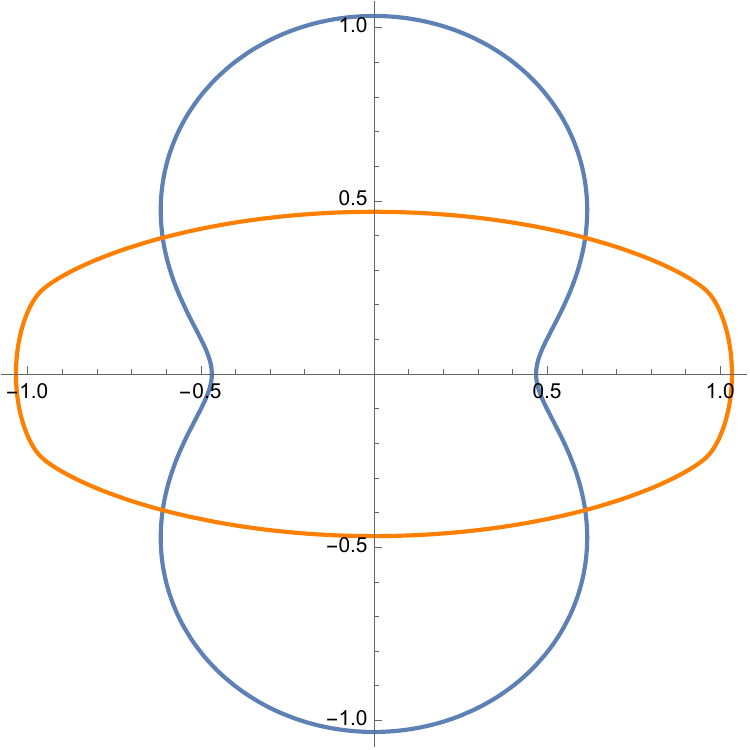}
    \end{minipage}
     \begin{minipage}[b]{.5\linewidth}
     \centering
    \includegraphics[scale= 0.45]{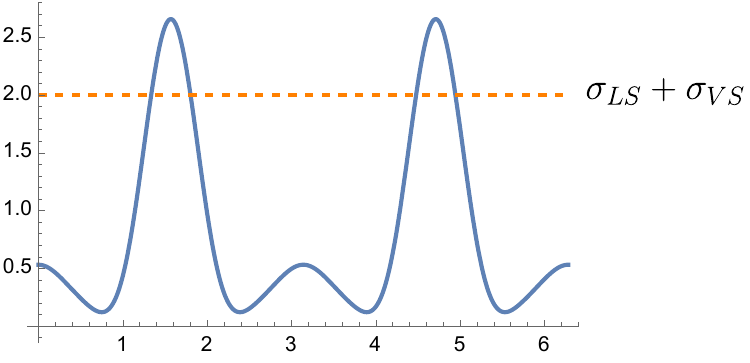}
    \end{minipage}
     \caption{Left: the blue line is the surface tension in polar coordinates as defined in (\ref{eq:sf}), and the orange line is the associated Wulff shape. Right: the violation of strong triangle inequality. The blue line is $\sigma_{VL} + \sigma''_{VL}$, and the orange dashed line is $\sigma_{LS} + \sigma''_{LS} + \sigma_{VS} + \sigma''_{VS}$, which is simply $\sigma_{LS} + \sigma_{VS}$.}
     \label{fig:surface_tension_wulff_prescribed_sf_mobility_examples}
\end{figure}

In this section, we investigate the convergence of the median filter scheme with both the surface tension $\sigma_{VL}$ {\em and} the mobility $\textcolor{black}{m}_{VL}$ of the vapor-liquid interface \textcolor{black}{as} prescribed anisotropic functions of inclination; this necessitates kernels concentrated on (at least) two circles.
Our choice of surface tensions of course respect\textcolor{black}{s} (\ref{eq:triangle}), but they violate (\ref{eq:st}) (\textcolor{black}{see Figure \ref{fig:surface_tension_wulff_prescribed_sf_mobility_examples}}) so that hypothesis of Claim \ref{claim:max} is not satisfied by our kernel construction (nevertheless convergence is observed):
\begin{equation}
\label{eq:sf}
    \sigma_{VL}(\theta) := \frac{7420\cos (2\theta) - 1316\cos(4\theta) + 180\cos(6\theta) - 25\cos(8\theta) + 21525}{26880}
\end{equation}
and the other two surface tensions are $\sigma_{LS} = \sigma_{VS} = 1$.\\

\begin{figure}
 \begin{minipage}[b]{.45\linewidth}
    \centering
    \includegraphics[scale= 0.4]{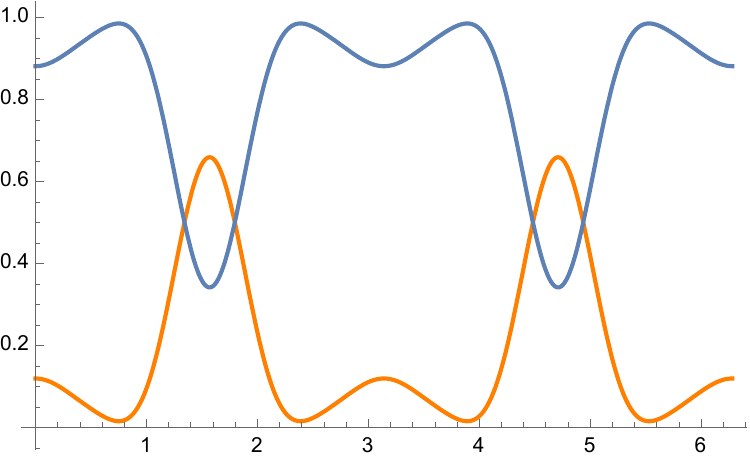}
    \end{minipage}
     \begin{minipage}[b]{.45\linewidth}
     \centering
    \includegraphics[scale= 0.4]{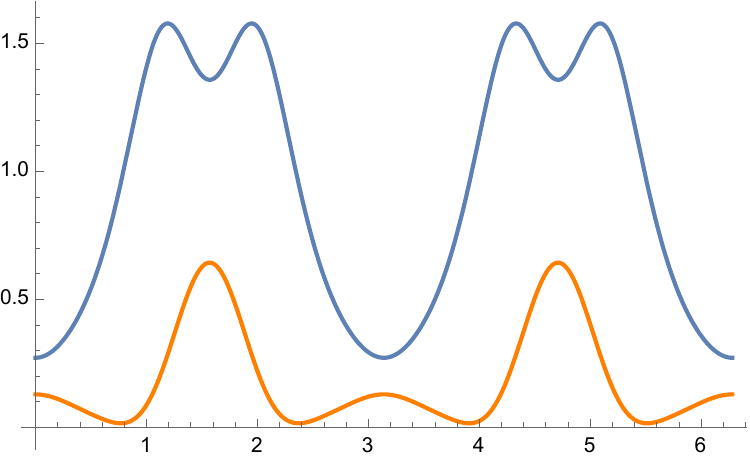}
    \end{minipage}
     \caption{Weight functions on the support of the kernels. The support of the kernels is chosen to be $\textcolor{black}{\partial B}_{\frac{1}{4}}(0)$ and $\textcolor{black}{\partial B}_{2}(0)$. The blue lines are the weight functions defined on $\textcolor{black}{\partial B}_{\frac{1}{4}}(0)$, and the orange lines are the weight functions defined on $\textcolor{black}{\partial B}_{2}(0)$. The x-axis for each figure represents the angle in radians. Left: mobility $\textcolor{black}{m(\theta)} = 1$. Right: mobility $\textcolor{black}{m}(\theta) = \frac{1}{2}+2\cos^4\theta$.}
        \label{fig:prescribed_sf_mobility_examples}
\end{figure}
\begin{figure}
 \begin{minipage}[b]{.45\linewidth}
    \centering
    \includegraphics[scale= 0.45]{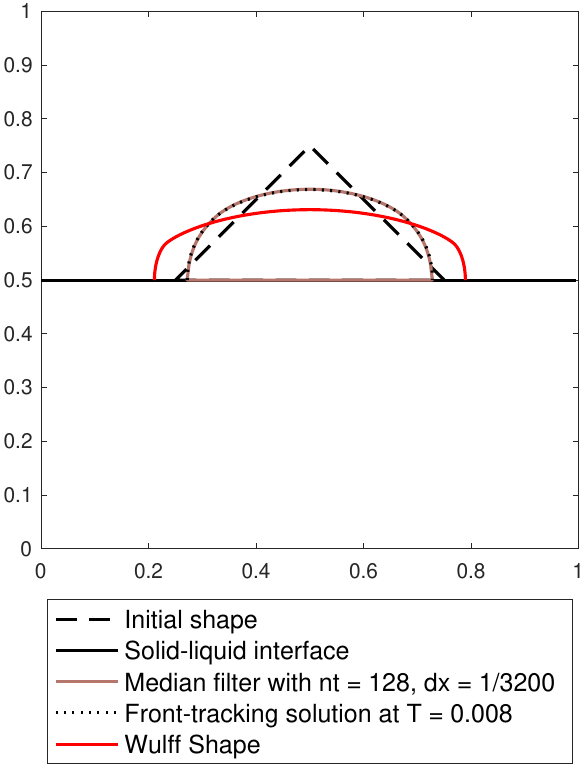}
    \end{minipage}
     \begin{minipage}[b]{.45\linewidth}
     \centering
    \includegraphics[scale= 0.45]{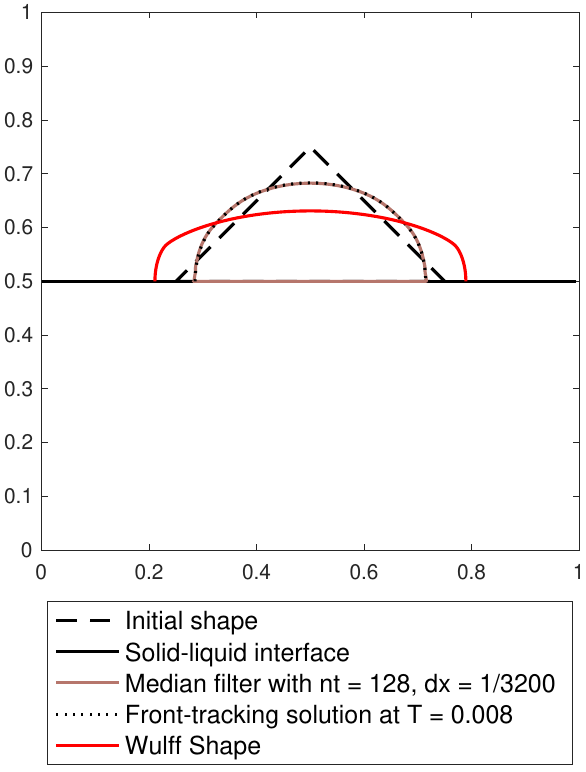}
    \end{minipage}
     \caption{Dynamics of droplets with the same surface tension $\sigma_{VL}$ but different mobilities, $\textcolor{black}{m}_1$ and $\textcolor{black}{m}_2$. The stationary shapes for both cases are the same, as they depend only on the prescribed surface tension $\sigma_{VL}$, while there are differences in their dynamics. Final time $T = 0.008$. Left: mobility $\textcolor{black}{m}_1 = 1$. Right: mobility $\textcolor{black}{m}_2(\theta) = \frac{1}{2} + 2\cos^4\theta$.}
        \label{fig:prescribed_sf_mobility_weight}
        
\end{figure}

The parametric expressions of the Wulff shape associated with this choice of $\sigma_{VL}$ is
\begin{equation}
      \begin{dcases}
        x(\theta) = \frac{65310\sin\theta + 14000\sin 3\theta + 5208\sin5\theta + 1125\sin7\theta + 175\sin9\theta}{53760}\\
        y(\theta) = \frac{20790\cos\theta + 840\cos3\theta + 2688\cos5\theta + 675\cos7\theta + 175\cos9\theta}{53760}
    \end{dcases}
\end{equation}

\begin{figure}
  \begin{minipage}[b]{.45\linewidth}
    \centering
    \begin{adjustbox}{width=\columnwidth,center}  
     \begin{tabular}[b]{c|c|c|c} 
     \hline
   \rule{0pt}{25pt} Number of time steps & $1/dx$ & $\textcolor{black}{L^{1}}$ Error & Order \\
    \hline
    \rule{0pt}{25pt} 2 & 50 & 0.010365   & -\\
     \hline
    \rule{0pt}{25pt}4 & 100 & 0.0062055   & 0.74011    \\
      \hline
   \rule{0pt}{25pt} 8 & 200 & 0.0033672     & 0.88202 \\
    \hline
   \rule{0pt}{25pt} 16 & 400 & 0.0017817  & 0.91832   \\
    \hline
   \rule{0pt}{25pt} 32 & 800 & 0.0010523     & 0.75968    \\
    \hline
    \rule{0pt}{25pt}64 & 1600 & 0.00077358   & 0.44392   \\
    \hline
   \rule{0pt}{25pt} 128  & 3200 & 0.00057721  & 0.42245   \\        \hline     
    \end{tabular}
    \end{adjustbox}
    \captionof{table}{Error Table for the case with mobility $\textcolor{black}{m}(\theta) = 1$.}
     \label{table:error_fitting_case1}
  \end{minipage}\hfill
  \begin{minipage}[b]{.45\linewidth}
  \centering
    \includegraphics[width=1.1\textwidth, center]{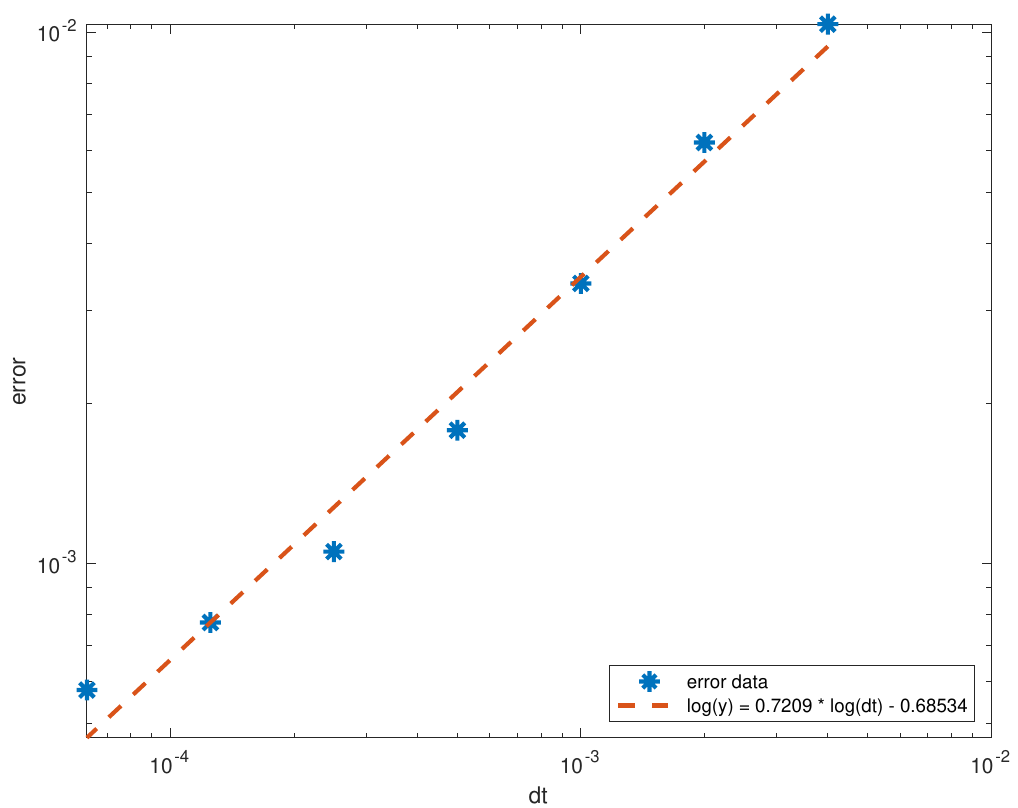}
    \captionof{figure}{Error plot with fitting line for the case with $\textcolor{black}{m}(\theta) = 1$. The slope of the log-log fitting line is 0.7209. }
     \label{fig:error_fitting_case1}
  \end{minipage}
  \label{fig:prescribed_sf_mobility_example1_error}
\end{figure}

\begin{figure}
  \begin{minipage}[b]{.45\linewidth}
    \centering
    \begin{adjustbox}{width=\columnwidth,center}  
     \begin{tabular}[b]{c|c|c|c} 
     \hline
   \rule{0pt}{25pt} Number of time steps & $1/dx$ & $\textcolor{black}{L^{1}}$ Error &  Order \\
    \hline
    \rule{0pt}{25pt} 2 & 50 & 0.014265   & -\\
     \hline
    \rule{0pt}{25pt}4 & 100 & 0.0089216   & 0.67711     \\
      \hline
   \rule{0pt}{25pt} 8 & 200 & 0.0046149     & 0.95099   \\
    \hline
   \rule{0pt}{25pt} 16 & 400 & 0.0021612   & 1.0945   \\
    \hline
   \rule{0pt}{25pt} 32 & 800 & 0.0011955    &  0.8542    \\
    \hline
    \rule{0pt}{25pt}64 & 1600 & 0.00066097   & 0.85495    \\
    \hline
   \rule{0pt}{25pt} 128  & 3200 & 0.00035217  & 0.90832  \\        \hline     
    \end{tabular}
    \end{adjustbox}
    \captionof{table}{Error Table for the case with mobility $\textcolor{black}{m}(\theta) = \frac{1}{2}+2\cos^4\theta$.}
    \label{table:error_fitting_case2}
  \end{minipage}\hfill
  \begin{minipage}[b]{.45\linewidth}
  \centering
    \includegraphics[width=1\textwidth, center]{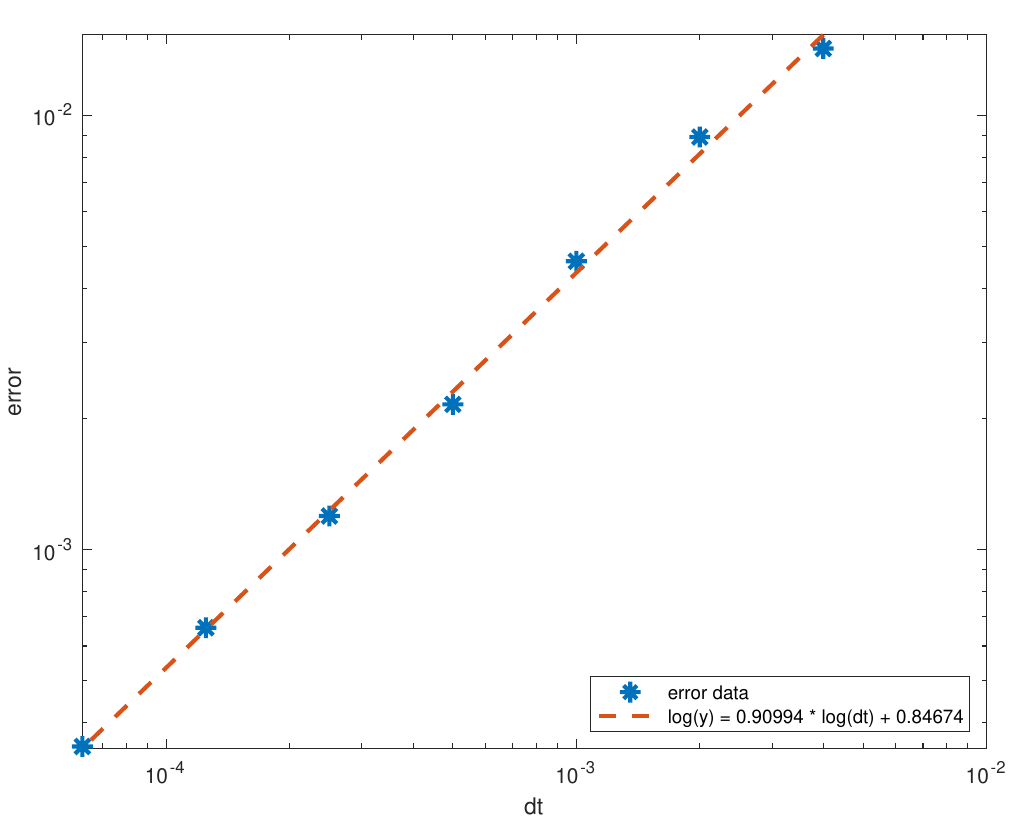}
    \captionof{figure}{Error plot with fitting line for the case with $\textcolor{black}{m}(\theta) = \frac{1}{2}+2\cos^4\theta$. The slope of the log-log fitting line is  0.9099. }
    \label{fig:error_fitting_case2}
  \end{minipage}
  \label{fig:prescribed_sf_mobility_example2_error}
\end{figure}

\noindent We present experiments with two different choices $\textcolor{black}{m}_1$ and $\textcolor{black}{m}_2$ for the mobility $\textcolor{black}{m}_{VL}$ of the only moving interface (namely, the vapor-liquid interface); they are given by:
\begin{equation}
\begin{dcases}
        \textcolor{black}{m}_1(\theta) := 1\\
        \textcolor{black}{m}_2(\theta) := \frac{1}{2}+2\cos^4\theta
\end{dcases}
\end{equation}
The two cases are expected to have the same stationary shape but different dynamics due to the difference between the mobilities. \textcolor{black}{Figure \ref{fig:prescribed_sf_mobility_examples} illustrates the different weight functions defined on the support of the kernels.}
We choose $T = 0.008$ to compare against the benchmark front-tracking solution, at which time the shape of the droplet is substantially different from the initial shape, but also far from the stationary (Wulff) shape, to test convergence of the dynamics and not just stationary states.
See Tables \ref{table:error_fitting_case1} and \ref{table:error_fitting_case2}, Figures \ref{fig:error_fitting_case1} and \ref{fig:error_fitting_case2} for details. 
Also see Figures \textcolor{black}{\ref{fig:prescribed_sf_mobility_weight}} for comparison of the best results in the two cases together with their Wulff shape and very accurate front-tracking solutions. 

\subsection{Fully anisotropic case}
\label{sec:numerical_fully_anisotropic}
In this section, we study the convergence of the vectorial median filter in the fully anisotropic setting, where the anisotropic surface tensions $\sigma_{VL}$, $\sigma_{VS}$, $\sigma_{LS}$, and mobility $\mu_{VL}$ are prescribed.
To this end, we select two pairs of surface tension and mobility—each pair being either normal-dependent or non-trivial (i.e., $\mu(\sigma+\sigma'') \neq \text{constant}$).
Note that the mobilities for the interfaces $\Gamma_{VS}$ and $\Gamma_{LS}$ are expected to have only minor effects on the numerical results, which diminish as resolution increases, since the interfaces $\Gamma_{VS}$ and $\Gamma_{LS}$ are stationary; any minor differences arise from various kernels resulting from the construction in Section~\ref{sec:Conv_kernel_construction}.
To test this hypothesis, we conducted experiments with two different mobilities, $\textcolor{black}{m}_{VS1}$ and $\textcolor{black}{m}_{VS2}$, for the vapor-solid interface $\Gamma_{VS}$.

\begin{figure}[h]
    \centering
\includegraphics[scale = 0.65]{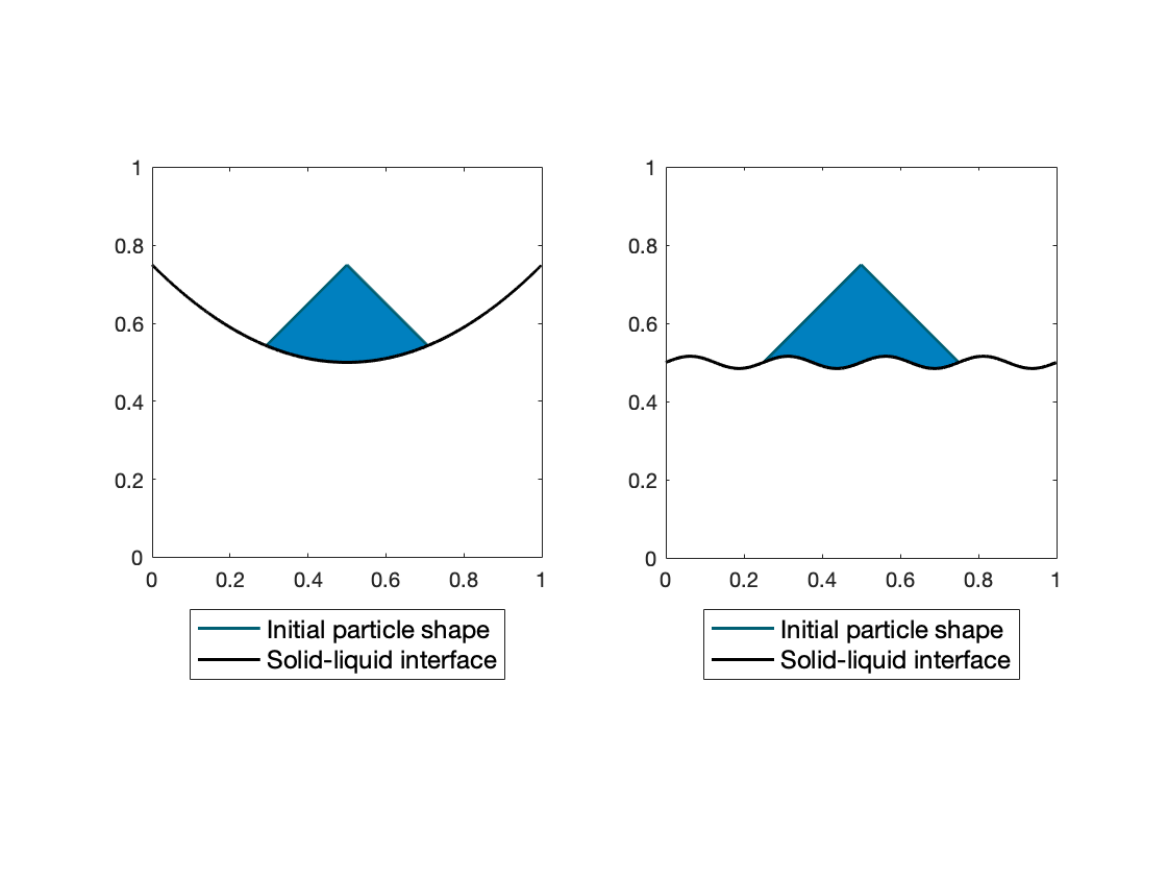}
    \caption{Initial configurations of different solid surfaces. Left: $y = (x-0.5)^2+0.5$. Right: $y = \frac{1}{64}\sin\big(8\pi(x-0.5)\big)+0.5$.}
    \label{fig:full_anisotropic_init}
\end{figure}

We consider the solid substrate to have two different \textcolor{black}{non-flat} shapes, as depicted in Fig.~\ref{fig:full_anisotropic_init}, which also displays the initial shape of the vapor-liquid interface $\Gamma_{VL}$.
The surface tension and mobility pairs for the three interfaces are as follows:

\begin{equation}
\label{eq:sf_mb_pair}
    \begin{split}
         \sigma_{VL}(\theta) = \sqrt{1+\sin^2(\theta+\frac{\pi}{4})}, \text{~~} &\textcolor{black}{m}_{VL}(\theta) = \sqrt{1+\cos^2\theta}\\
         \sigma_{LS}(\theta) = \sqrt{1+\sin^2(\theta+\frac{\pi}{3})}, \text{~~} &\textcolor{black}{m}_{LS}(\theta) = \frac{(3 + \sin(\frac{\pi}{6} + 2\theta))^{\frac{3}{2}}}{4\sqrt{2}}\\
         \sigma_{VS}(\theta) = \sqrt{1+\sin^2(\theta+\frac{\pi}{8})},\text{~~}&\textcolor{black}{m}_{VS1}(\theta) = \frac{(3 - \cos(\frac{\pi}{4} + 2\theta))^{\frac{3}{2}}}{4\sqrt{2}} \\&\textcolor{black}{m}_{VS2}(\theta) = \sqrt{1+\cos^2(\theta+\frac{\pi}{8})}
    \end{split}
\end{equation}

\begin{figure}[h]

    \centering
    \includegraphics[scale = 0.45]{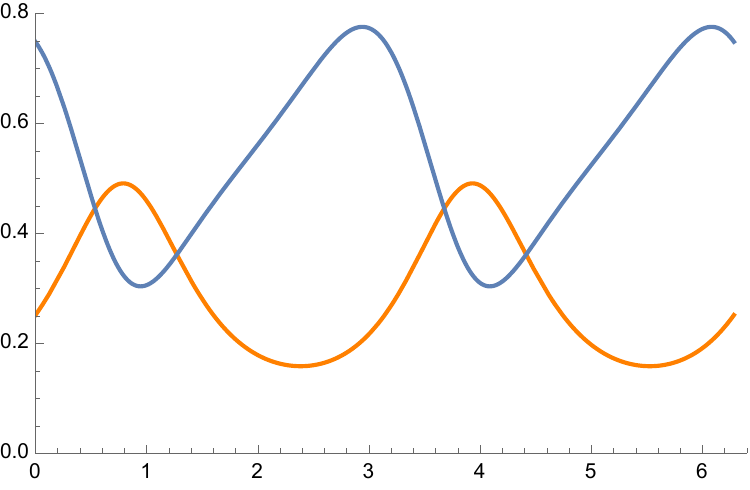}
    \hspace{1cm}
    \includegraphics[scale = 0.45]{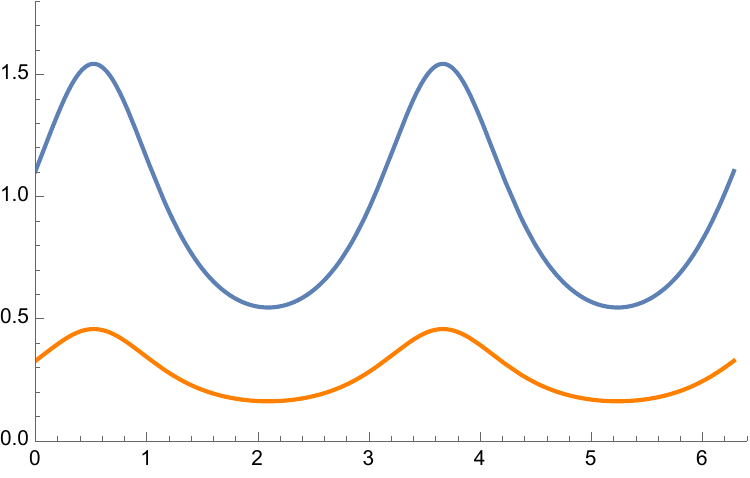}
  
    \includegraphics[scale = 0.45]
    {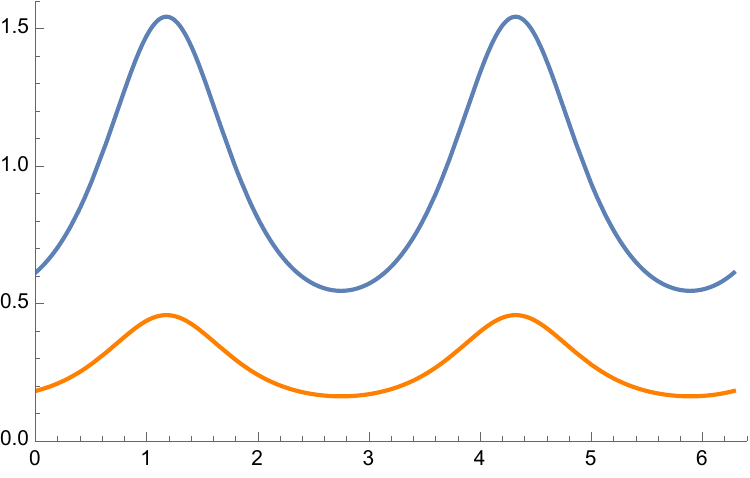}
     \hspace{1cm}
    \includegraphics[scale = 0.45]{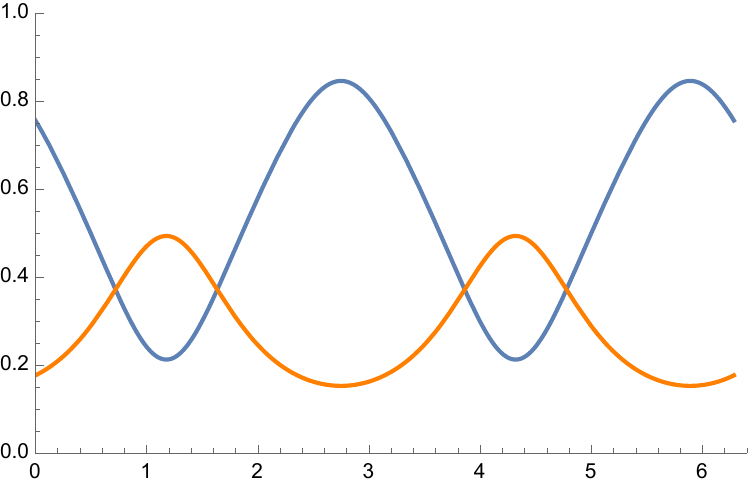}
    \caption{Weight functions of kernels in $\theta$-coordinates associated with different pairs of surface tension and mobility. In all the plots, the blue lines represent the weight functions on the circle with a radius of $\frac{1}{3}$, while the orange lines correspond to those on the circle with a radius of $2$. The plots in the first row depict the weight functions for the surface tension-mobility pairs $(\sigma_{VL}, \textcolor{black}{m}_{VL})$ and $(\sigma_{LS}, \textcolor{black}{m}_{LS})$, respectively. The two plots on the bottom row represent the weight functions for the surface tension $\sigma_{VS}$ with mobilities $\textcolor{black}{m}_{VS1}$ and $\textcolor{black}{m}_{VS2}$, respectively, as defined in Eq.~(\ref{eq:sf_mb_pair}).}
    \label{fig:fully_anisotropic_weights}.
\end{figure}

 In the kernel construction (see Section~\ref{sec:Conv_kernel_construction}), we have chosen the radii of the two circles to be $R_1 = \frac{1}{3}$ and $R_2 = 2$.
Polar plots of the corresponding weight functions for the kernel constructed for each pair of surface tension and mobility are shown in Figure~\ref{fig:fully_anisotropic_weights}.

\begin{figure}[h]
    \centering
    \includegraphics[scale = 0.35]{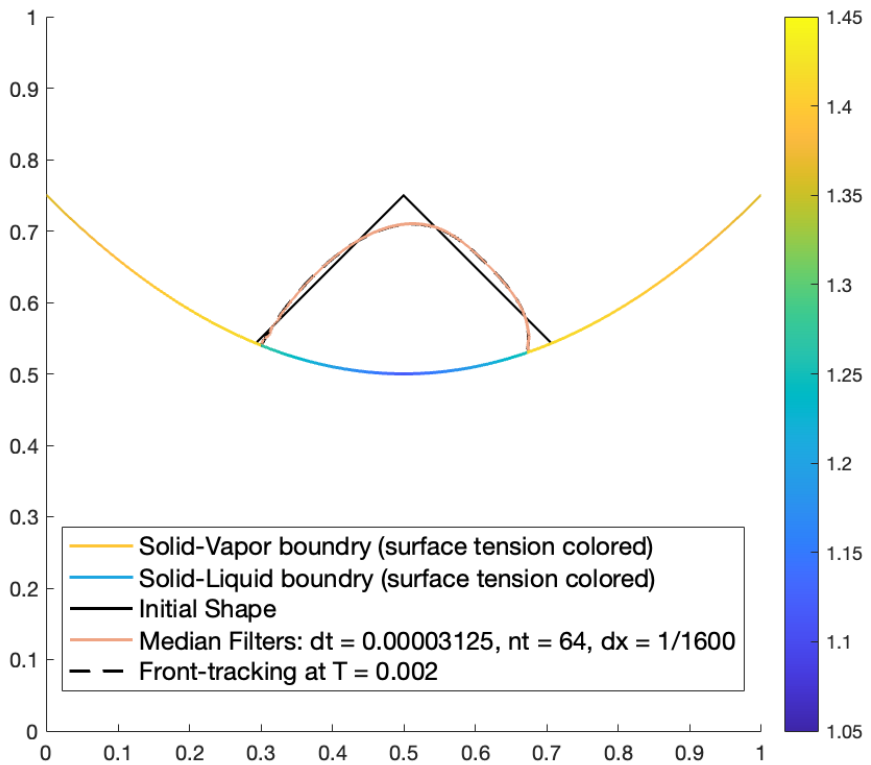}
    \hspace{1cm}
     \includegraphics[scale = 0.35]{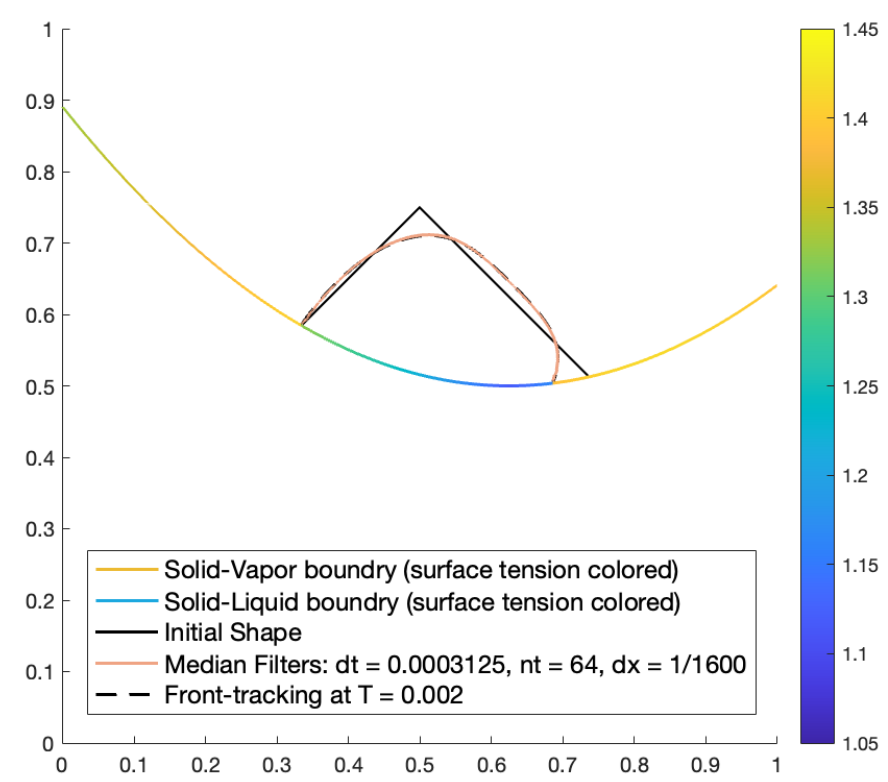}
    \caption{Median filter solution for the first fully anisotropic case with \textcolor{black}{substrate surface given by different parabola}, compared with the front-tracking solution. Here, the color bar indicates the surface tensions on the solid surface, the black solid line represents the initial shape of the droplet, the dashed line is the front-tracking solution, and the pink line is the result obtained by the median filter scheme, as shown in Figure~\ref{fig:fully_anisotropic_weights}.}
    \label{fig:fully_anisotropic_squarebry}
\end{figure}

\begin{figure}[h]
    \centering
    \includegraphics[scale = 0.35]{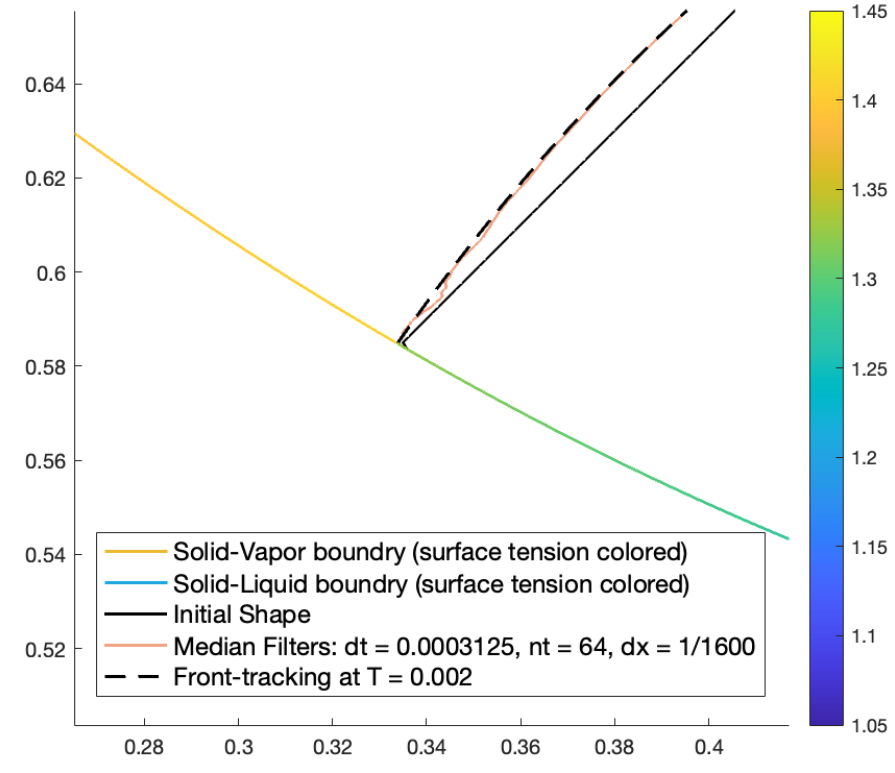}
    \hspace{1cm}
     \includegraphics[scale = 0.35]{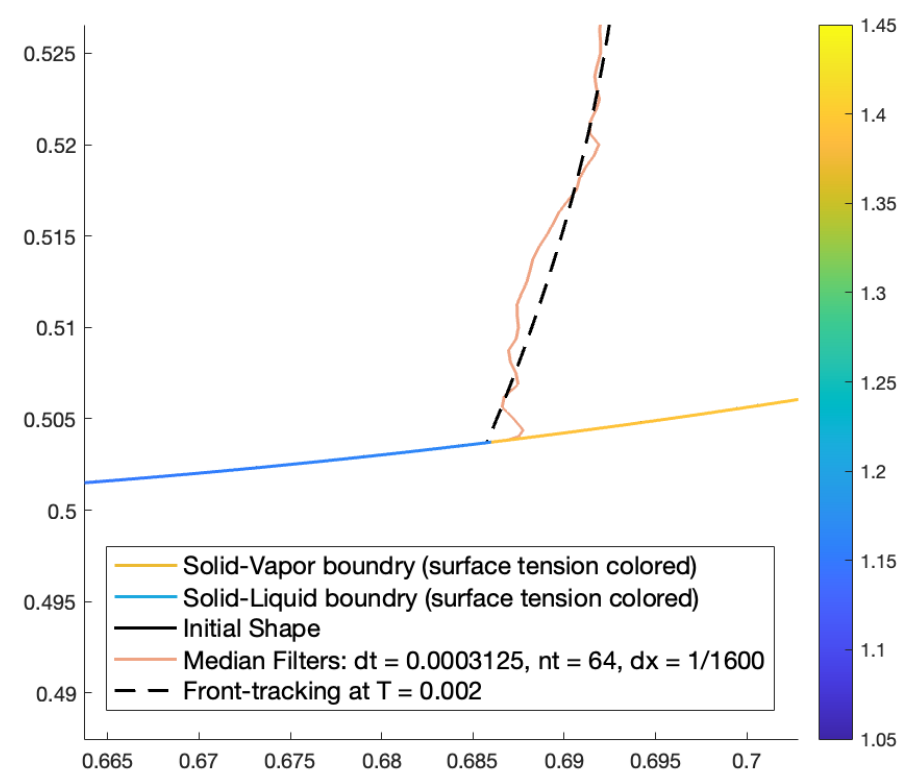}
    \caption{\textcolor{black}{Parabolic substrate surface}: details of contact angles}
    \label{fig:fully_anisotropic_squarebry_angles}
\end{figure}

\begin{figure}
  \begin{minipage}[b]{.45\linewidth}
    \centering
    \begin{adjustbox}{width=\columnwidth,center}  
     \begin{tabular}[b]{c|c|c|c} 
     \hline
   \rule{0pt}{30pt} Number of time steps & $1/dx$ & $\textcolor{black}{L^{1}}$ Error &  Order \\
    \hline
    \rule{0pt}{30pt}2 & 50 &  0.0058955   & -\\
     \hline
    \rule{0pt}{30pt}4 & 100 &  0.0032622 & 0.85    \\
      \hline
   \rule{0pt}{30pt}  8 & 200 &  0.0018778    & 0.80 \\
    \hline
   \rule{0pt}{30pt} 16 & 400 &  0.00093219 & 1.01   \\
    \hline
   \rule{0pt}{30pt} 32 & 800 &0.00052523     & 0.83   \\
    \hline
    \rule{0pt}{30pt} 64 & 1600 &  0.0003411   & 0.62   \\
    \hline     
    \end{tabular}
    \end{adjustbox}
    \captionof{table}{Error table for the initial configuration with a symmetric quadratic solid surface with the first solid-vapor mobility $\textcolor{black}{m}_{SV1}$ at $T = 0.002$}
    \label{table:error_table_for_sym_quad}
  \end{minipage}\hfill
  \begin{minipage}[b]{.45\linewidth}
  \centering
    \includegraphics[width=1.1\textwidth, center]{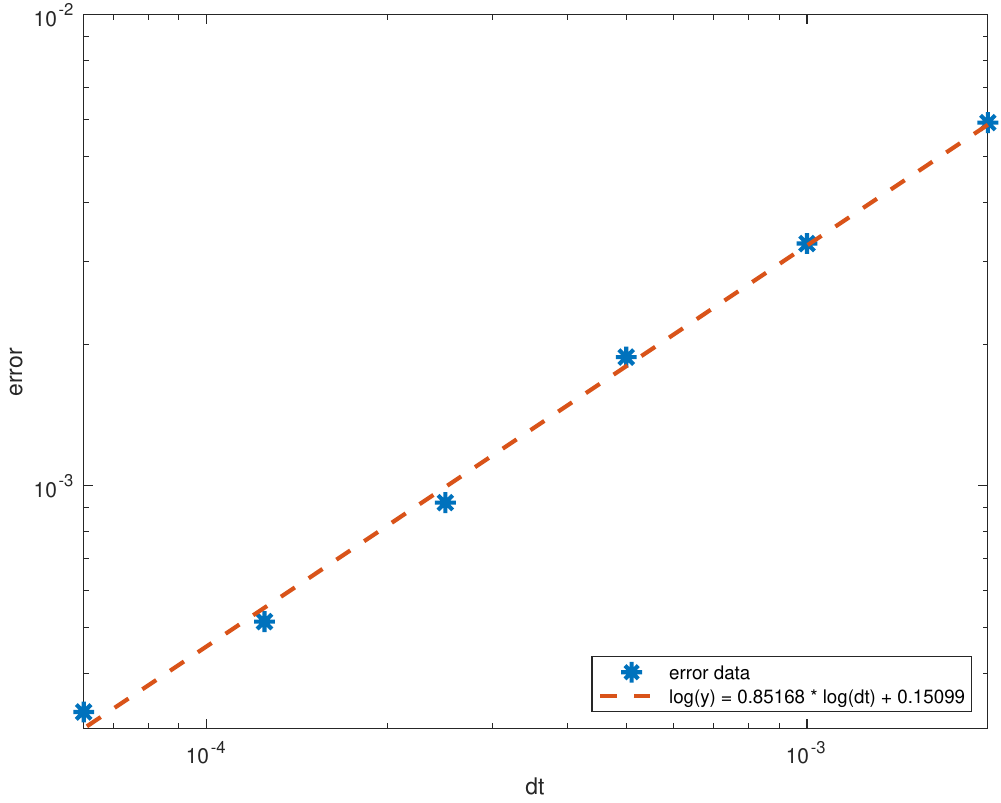}
    \captionof{figure}{Error plot with fitting line for the case with quadratic solid-vapor/liquid interface. The slope of the log-log fitting line is 0.85168. }
    \label{fig:error_table_for_sym_quad}
  \end{minipage}
  \label{fig:fully_anisotropic_xsquare_sym}
\end{figure}

\begin{figure}
  \begin{minipage}[b]{.45\linewidth}
    \centering
    \begin{adjustbox}{width=\columnwidth,center}  
     \begin{tabular}[b]{c|c|c|c} 
     \hline
   \rule{0pt}{30pt} Number of time steps & $1/dx$ & $\textcolor{black}{L^{1}}$ Error & Order \\
    \hline
    \rule{0pt}{30pt}2 & 50 &  0.0063712   & -\\
     \hline
    \rule{0pt}{30pt}4 & 100 &  0.0035687 & 0.84    \\
      \hline
   \rule{0pt}{30pt}  8 & 200 &  0.0018259    & 0.97 \\
    \hline
   \rule{0pt}{30pt} 16 & 400 &  0.001162 & 0.65   \\
    \hline
   \rule{0pt}{30pt} 32 & 800 &0.00078184     & 0.57   \\
    \hline
    \rule{0pt}{30pt} 64 & 1600 &  0.00066438   & 0.23   \\
    \hline     
    \end{tabular}
    \end{adjustbox}
    \captionof{table}{Error table for the initial configuration with a non-symmetric quadratic solid surface with the first solid-vapor mobility $\textcolor{black}{m}_{SV1}$ at $T = 0.002$.}
    \label{table:error_table_for_nonsym_quad}
  \end{minipage}\hfill
  \begin{minipage}[b]{.45\linewidth}
  \centering
    \includegraphics[width=1.1\textwidth, center]{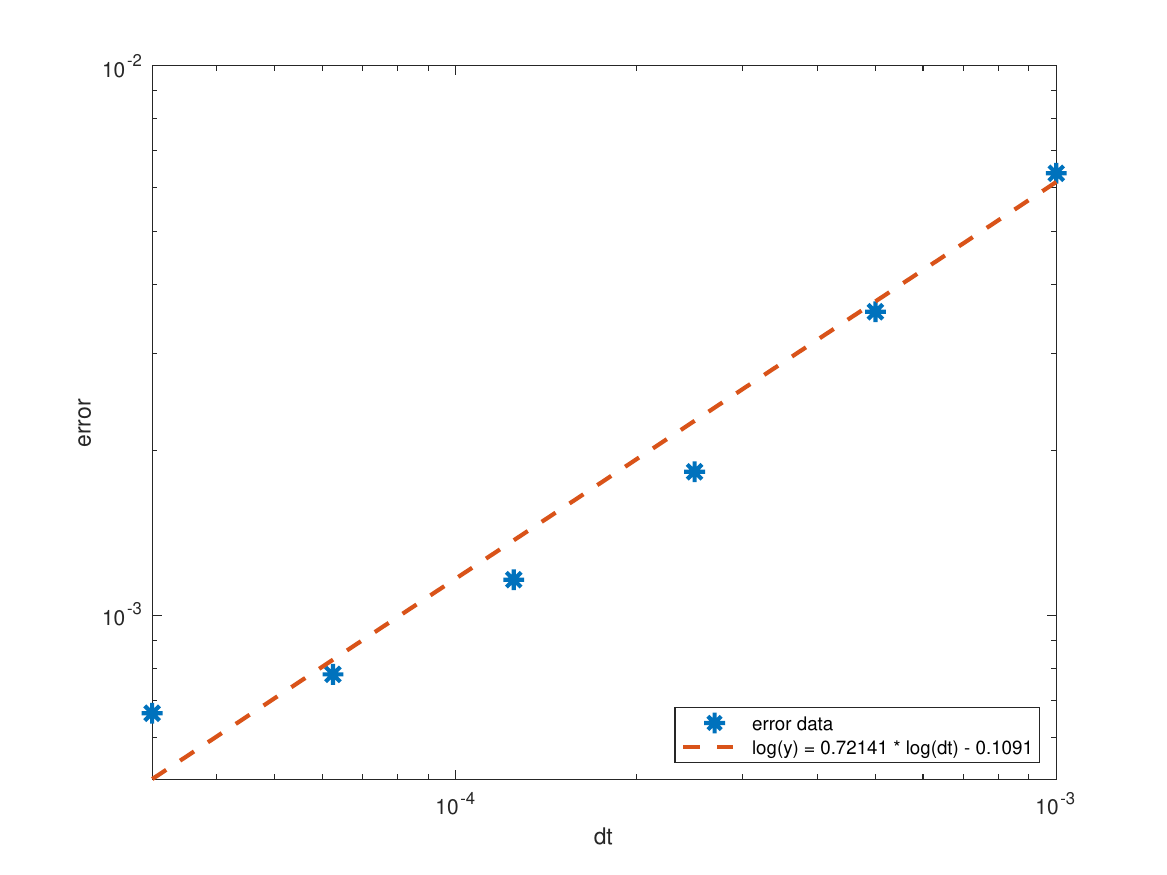}
    \captionof{figure}{Error plot with the fitting line for the case with quadratic solid-vapor/liquid interface. The slope of the log-log fitting line is 0.72141. }
    \label{fig:error_table_for_nonsym_quad}
  \end{minipage}
  \label{fig:fully_anisotropic_xsquare_nonsym}
\end{figure}

\textcolor{black}{For our first fully anisotropic test, we choose non-flat substrates given by a quadratic functions (parabolas) as show in the left panel of Figure \ref{fig:full_anisotropic_init}. Figure \ref{fig:fully_anisotropic_squarebry_angles} shows a zoom of the contact angles starting from the initial coidition shown in the right hand panel of Figure \ref{fig:fully_anisotropic_squarebry}}. We choose the mobility for the solid-vapor interface to be \textcolor{black}{$m_{VS1}$ from (\ref{eq:sf_mb_pair})} and compare it with the front-tracking solution at the final time $T = 0.002$ using a very fine $dt$ and $dx$. The visualization is shown in Figure~\ref{fig:fully_anisotropic_squarebry}, and the errors and order of errors (in time) are shown in Tables~\ref{table:error_table_for_sym_quad} and \ref{table:error_table_for_nonsym_quad}, as well as in Figures~\ref{fig:error_table_for_sym_quad} and \ref{fig:error_table_for_nonsym_quad}, respectively.

\begin{figure}
  \begin{minipage}[b]{.45\linewidth}
    \centering
    \begin{adjustbox}{width=\columnwidth,center}  
     \begin{tabular}[b]{c|c|c|c} 
     \hline
   \rule{0pt}{30pt} Number of time steps & $1/dx$ & $\textcolor{black}{L^{1}}$ Error & Order \\
    \hline
    \rule{0pt}{25pt}1 & 25 &  0.013391    & -\\
     \hline
    \rule{0pt}{25pt} 2 & 50 &  0.0094366    & 0.50    \\
      \hline
   \rule{0pt}{25pt}  4 & 100 &  0.0059839    &  0.66 \\
    \hline
   \rule{0pt}{25pt}  8 & 200 & 0.0030785  &  0.96  \\
    \hline
   \rule{0pt}{25pt}16 & 400 & 0.0016449   &  0.90    \\
    \hline
    \rule{0pt}{25pt} 32 & 800 & 0.0012686   & 0.37   \\
    \hline     
    \rule{0pt}{25pt} 64 & 1600 &  0.001099   & 0.21\\
    \hline
    \rule{0pt}{25pt} 128 & 3200 & 9.0307e-04 & 0.28
    \end{tabular}
    \end{adjustbox}
    \captionof{table}{Error table for the initial configuration with a sinusoidal solid surface with the first solid-vapor mobility $\textcolor{black}{m}_{SV1}$ at $T = 0.004$.}
     \label{table:fully_anisotropic_x_sin1}
  \end{minipage}\hfill
  \begin{minipage}[b]{.5\linewidth}
  \centering
    \includegraphics[width=1.1\textwidth, center]{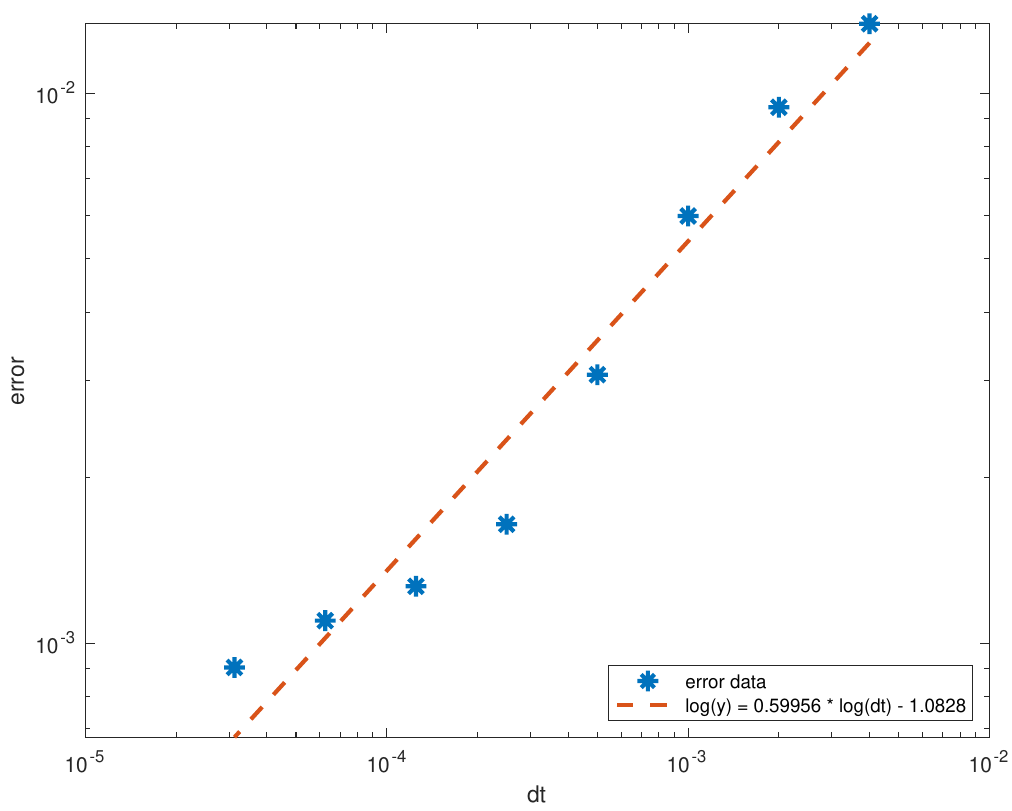}
    \captionof{figure}{Error plot with fitting line for the case with sinusoidal solid-vapor/liquid interface. The slope of the log-log fitting line is 0.59956. }
     \label{fig:fully_anisotropic_x_sin1}
  \end{minipage}
 
\end{figure}

\begin{figure}
  \begin{minipage}[b]{.45\linewidth}
    \centering
    \begin{adjustbox}{width=\columnwidth,center}  
     \begin{tabular}[b]{c|c|c|c} 
     \hline
   \rule{0pt}{30pt} Number of time steps & $1/dx$ & $\textcolor{black}{L^{1}}$ Error & Order \\
    \hline
    \rule{0pt}{25pt} 1 & 25 &  0.013062    & -\\
     \hline
    \rule{0pt}{25pt} 2 & 50 &  0.0095881    & 0.45   \\
      \hline
   \rule{0pt}{25pt}   4 & 100 & 0.0062323     &  0.62 \\
    \hline
   \rule{0pt}{25pt} 8 & 200 &  0.0032225   &  0.95  \\
    \hline
   \rule{0pt}{25pt} 16 & 400 & 0.0017067    &  0.92   \\
    \hline
    \rule{0pt}{25pt} 32 & 800 & 0.0012532  & 0.45   \\
    \hline     
    \rule{0pt}{25pt}   64 & 1600 &   0.0010780    & 0.22\\
    \hline
    \rule{0pt}{25pt}  128 & 3200 & 9.0306e-04 & 0.26
    \end{tabular}
    \end{adjustbox}
    \captionof{table}{Error table for the initial configuration with a sinusoidal solid surface with the second solid-vapor mobility $\textcolor{black}{m}_{SV2}$ at $T = 0.004$.}
    \label{table:fully_anisotropic_x_sin2}
  \end{minipage}\hfill
  \begin{minipage}[b]{.5\linewidth}
  \centering
    \includegraphics[width=1.1\textwidth, center]{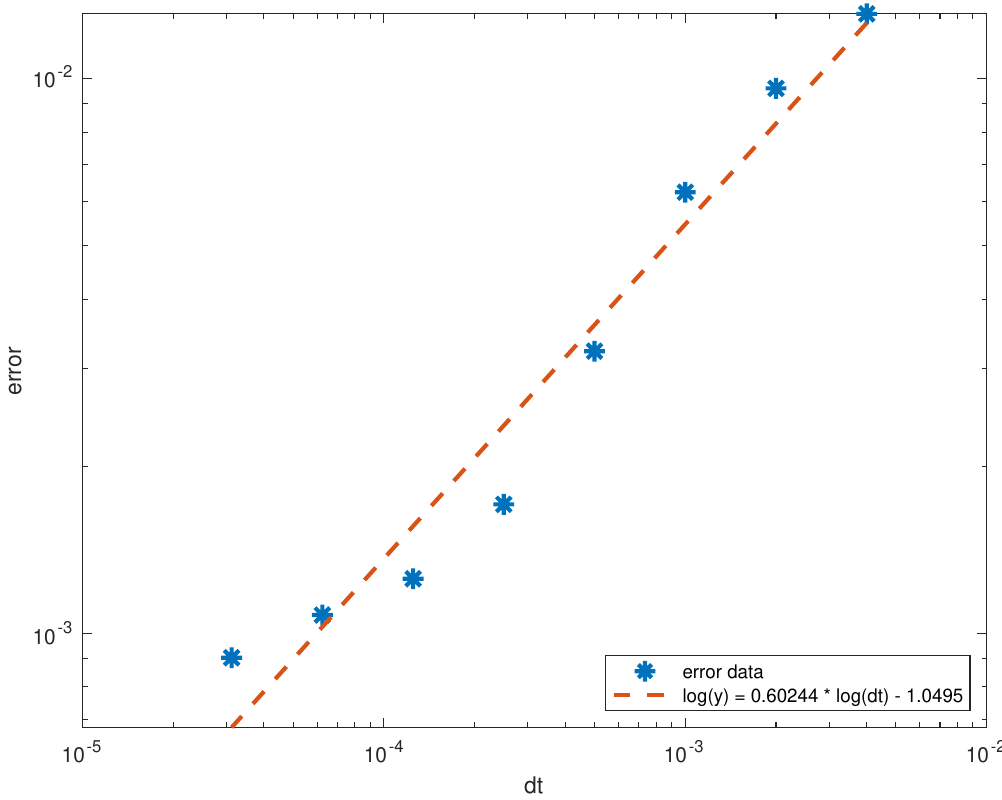}
    \captionof{figure}{Error plot with fitting line for the case with sinusoidal solid-vapor/liquid interface. The slope of the log-log fitting line is 0.60244. }
     \label{fig:fully_anisotropic_x_sin2}
  \end{minipage}
\end{figure}

\begin{figure}
    \centering
    \includegraphics[scale = 0.35]{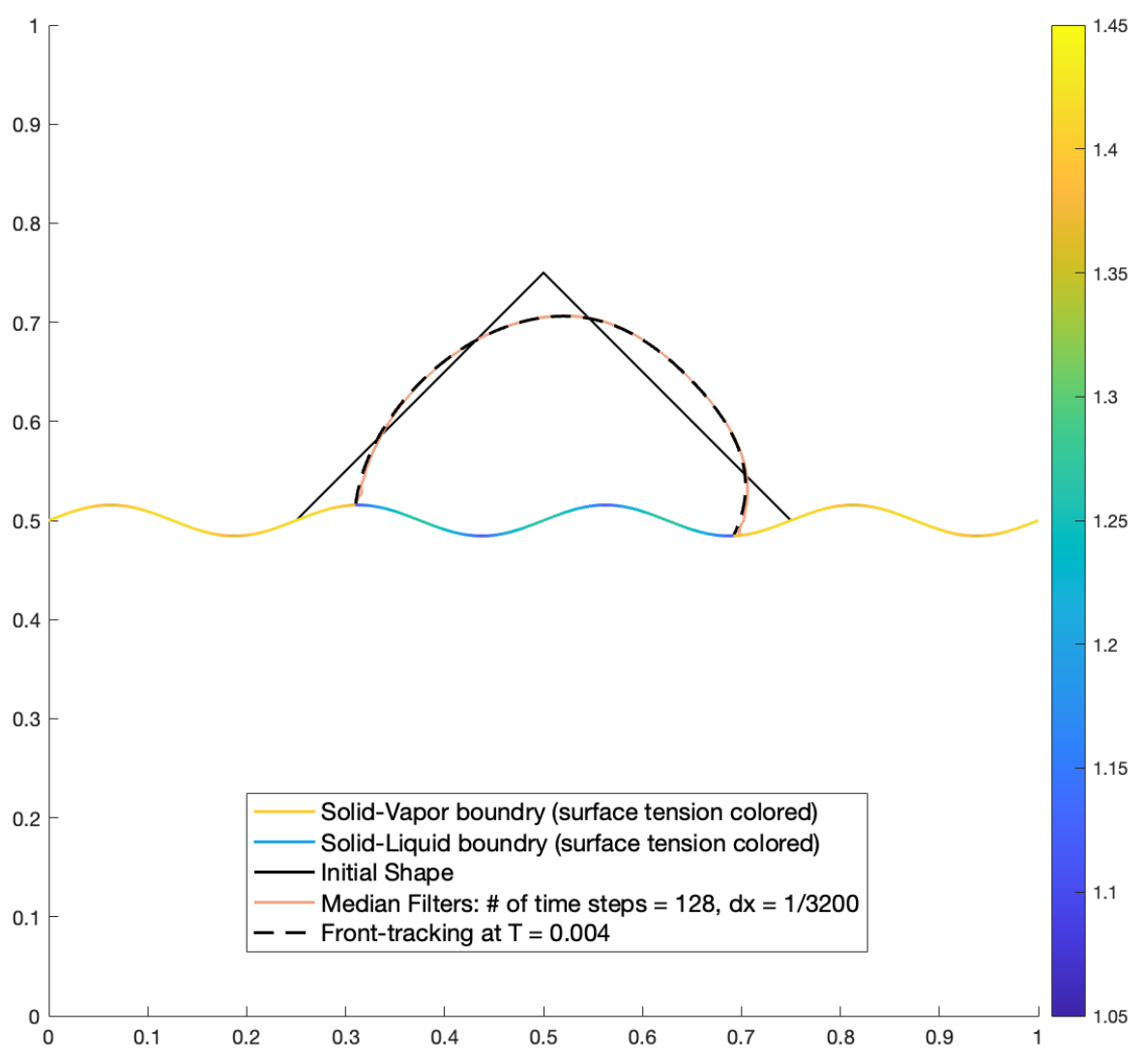}
    \caption{Fully anisotropic test with sinusoidal boundry}
    \label{fig:fully_anisotropic_x_sin_best}
\end{figure}

\begin{figure}
    \centering
    \includegraphics[scale = 0.25]{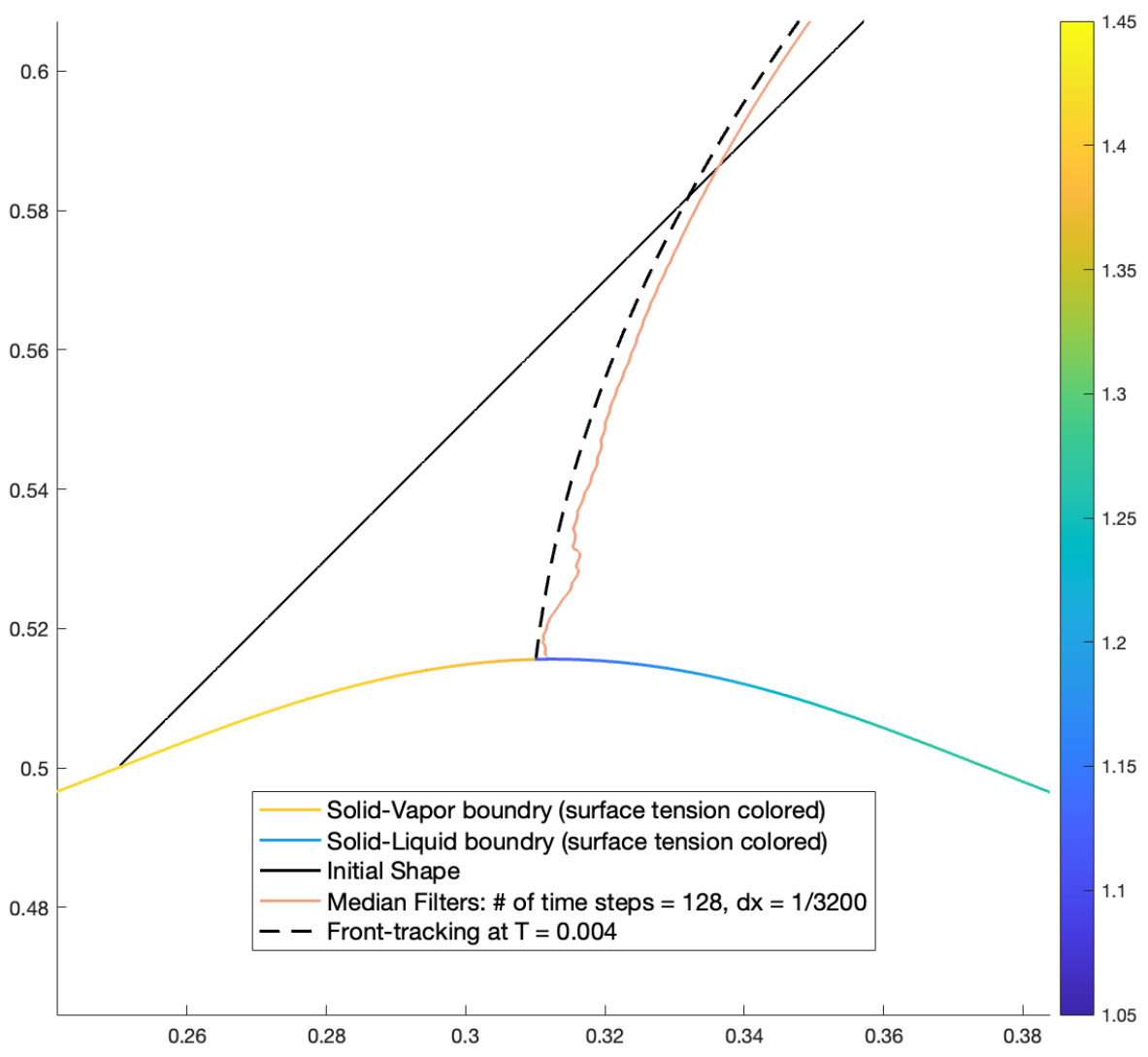}
    \hspace{1cm}
    \includegraphics[scale = 0.25]{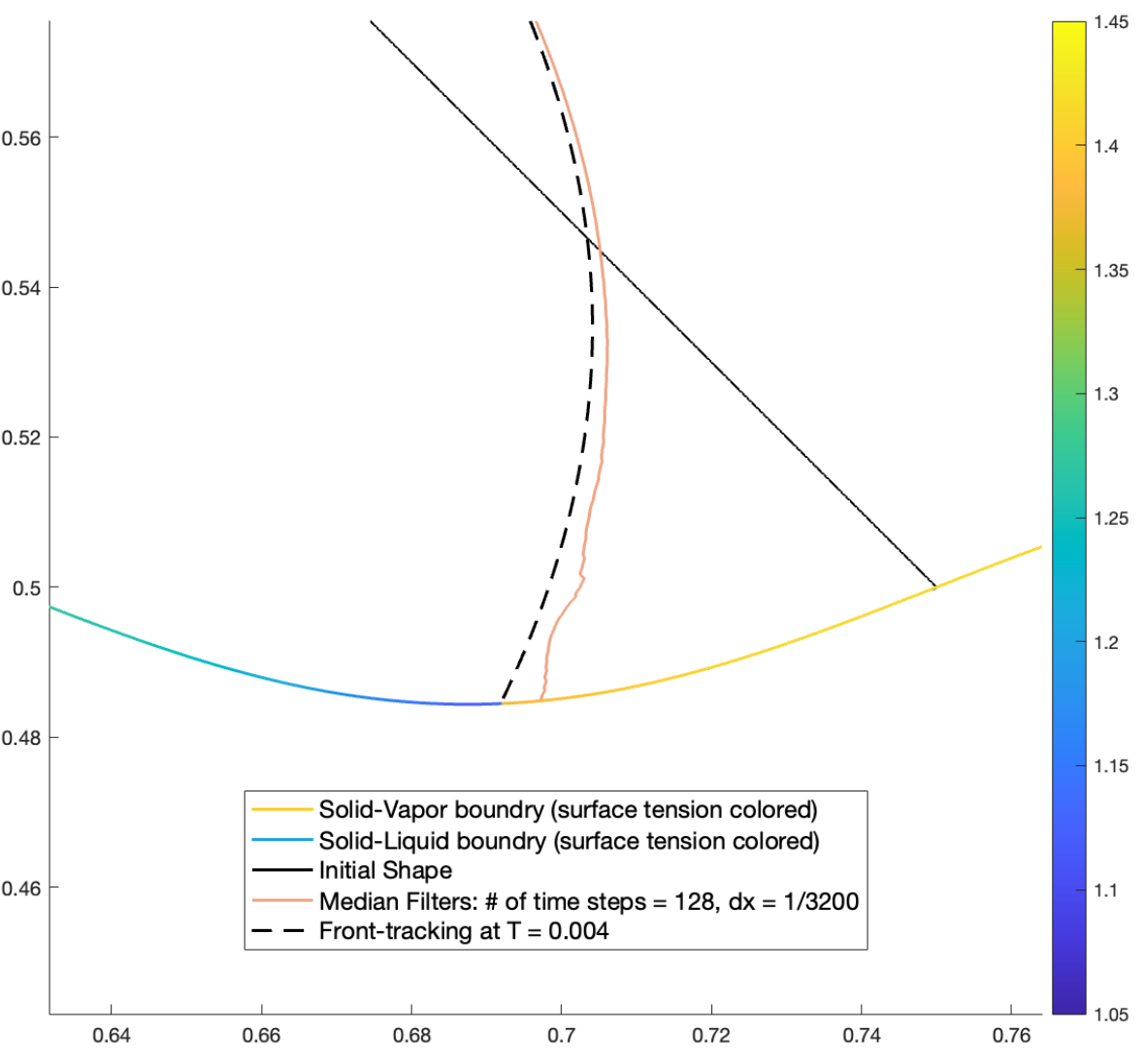}
    \caption{Fully anisotropic dynamic convergence test with sinusoidal boundary with contact angle details. There are still minor oscillations around the contact angles, but the error tables still show the desired convergence rate.}
    \label{fig:fully_anisotropic_x_sin_best_angles}
\end{figure}

For \textcolor{black}{our second test we use $y = \frac{1}{64}\sin(8\pi x)$ as the shape of the substrate, show in the right hand panel of Figure \ref{fig:full_anisotropic_init}}. We ran two simulations with different mobilities, $\textcolor{black}{m}_{VS1}$ and $\textcolor{black}{m}_{VS2}$, for the solid-vapor interface as defined in Eq.~(\ref{eq:sf_mb_pair}). The similarity in the errors reported in Tables \ref{table:fully_anisotropic_x_sin1} and \ref{table:fully_anisotropic_x_sin2}, Figures \ref{fig:fully_anisotropic_x_sin1} and \ref{fig:fully_anisotropic_x_sin2} for the two choices of mobilities substantiates the minor effects previously discussed at the beginning of this section. \textcolor{black}{See Figures \ref{fig:fully_anisotropic_x_sin_best} and \ref{fig:fully_anisotropic_x_sin_best_angles} for the best result we achieved, along with a comparison to the highly accurate contact angles obtained by front-tracking.}

\subsection{Topological Changes}
\label{sec:numerical_topological_changes}
\noindent One of the most notable advantages of the level set method is the ability to handle topological changes without ad-hoc numerical surgery, compared to the explicit methods. We performed two experiments involving topology change of the volume preserved liquid droplets evolving on a substrate with a sinusoidal solid surface. Both cases start with two separated droplets with an anisotropic surface tension on the interface with vapor:
\begin{equation}
\label{eq:anisotropy}
    \sigma_{VL}(\theta):= \sqrt{1+\sin^2(\theta-\frac{\pi}{3})}
\end{equation}
In the first case, the other two surface tensions are $\sigma_{LS} = 1, \sigma_{VS} = 1.5$, and in the second case, the other two surface tensions are $\sigma_{LS} = 1.5, \sigma_{VS} = 1$. The first one is expected to have wetting behavior while the second one is expected to have dewetting behavior. We used spatial grid size $dx = 0.0025$ and time step $\delta t = 0.0001$ for both cases. 
\textcolor{black}{We utilized narrow-banding to save running time, updating the level set functions only in a tubular neighborhood of their zero level sets.}\\

\begin{figure}
    \centering
    \includegraphics[scale = 0.6]{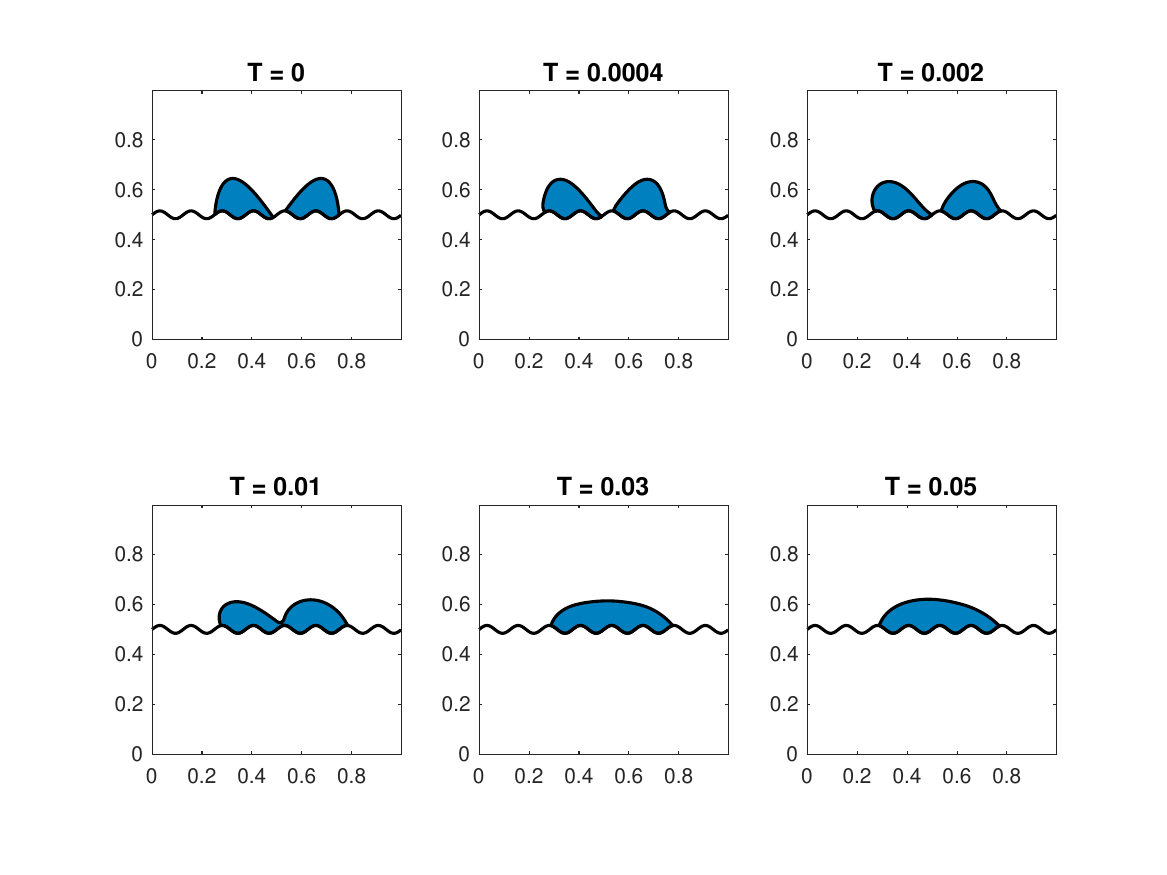}
    \caption{Wetting of separated droplets with anisotropic surface tension}
    \label{fig:two_comp_wetting}
\end{figure}

\begin{figure}
    \centering
    \includegraphics[scale = 0.6]{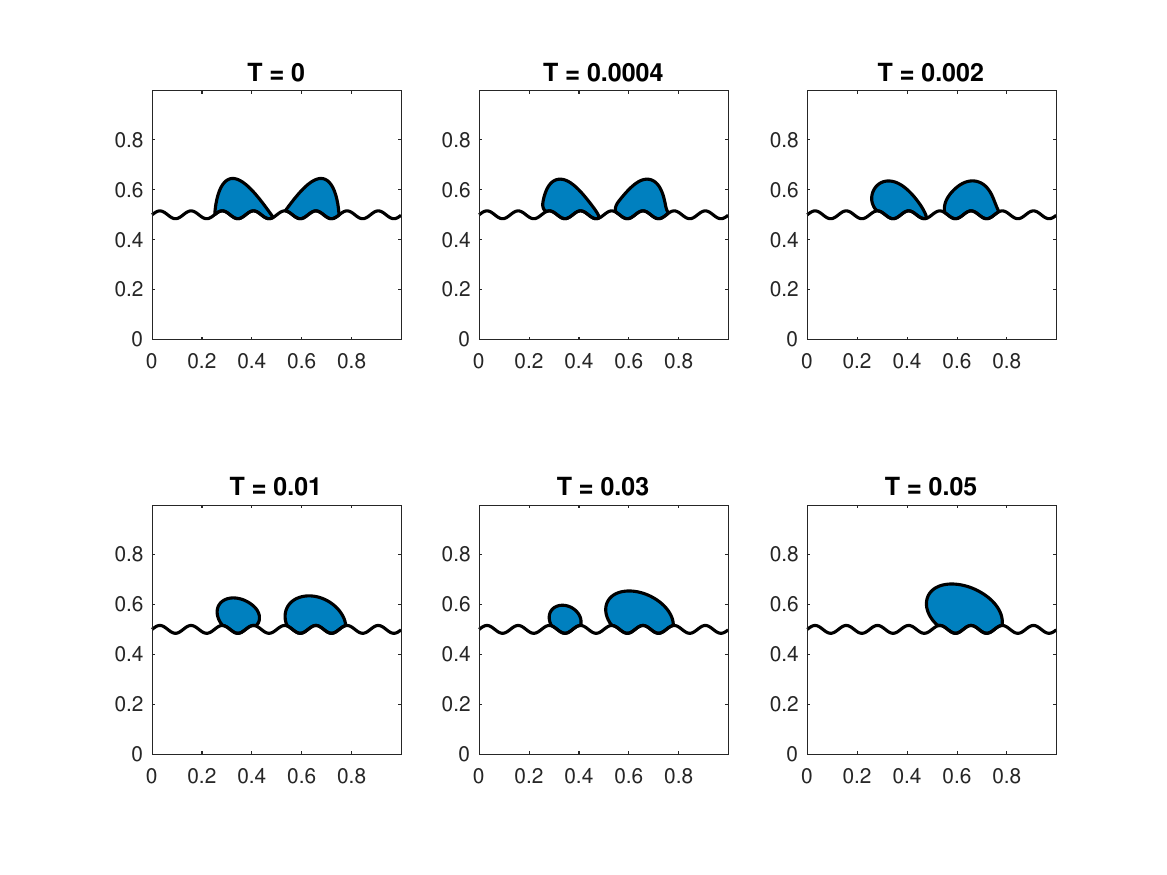}
    \caption{Dewetting of separated droplets with anisotropic surface tension.}
    \label{fig:two_comp_dewetting}
\end{figure}

\noindent As shown in Figure \ref{fig:two_comp_wetting}, the two droplets merge at around $T = 0.01$ and then evolve towards the Wulff shape. As for the dewetting case, the two droplets evolve separately until one of them disappears. See Figure \ref{fig:two_comp_dewetting} for the details.

\begin{figure}[h]
    \centering
    \includegraphics[scale=0.6]{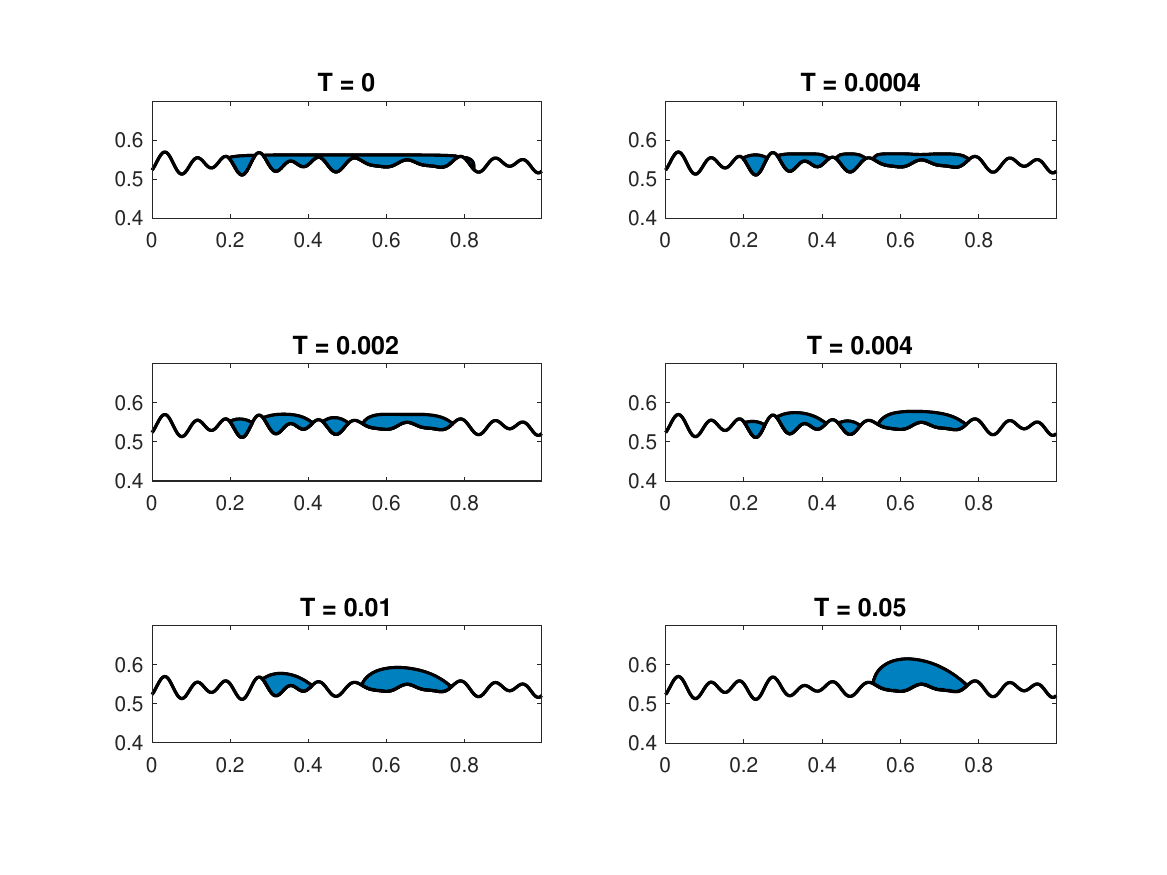}
    \caption{\textcolor{black}{Evolution of a long thin particle with anisotropic surface tension $\sqrt{1+\sin^2(\theta-\frac{\pi}{3})}$, $\sigma_{LS} = 1, \sigma_{VS} = 1$.}}
    \label{fig:thinfilm}
\end{figure}

\textcolor{black}{Finally, we test our algorithm on the evolution of an initially long and thin particle on a rough substrate with the same anisotropy $\sigma_{VL}(\theta)$ defined in (\ref{eq:anisotropy}) and $\sigma_{LS} = 1, \sigma_{VS} = 1$. Multiple topological changes occur during the evolution; see Figure \ref{fig:thinfilm}.}

\section{Conclusion}
Vectorial median filters are natural level set versions of threshold dynamics algorithms.
There have been many advances in our understanding of threshold dynamics in recent years, leading to its extension to particularly challenging problems such as anisotropic, multiphase motion by mean curvature of networks.
One drawback of threshold dynamics algorithms is their low accuracy when implemented naively on uniform spatial grids; this lack of accuracy, which is an unfortunate byproduct of representing interfaces via discontinuous (characteristic) functions, can lead to numerical artifacts such as pinning if spatial resolution is insufficient.
The vectorial median filter versions of these algorithms alleviate this difficulty by representing interfaces via regular level set functions that can be interpolated to locate the interface at subgrid accuracy, while maintaining consistency with threshold dynamics in principle and thereby benefiting from its existing theory and generality.
In this paper, we explored a median filter version of threshold dynamics in the context of wetting / dewetting problems with anisotropic surface tensions and mobilities.
It is a natural level set version of the anisotropic threshold dynamics from \cite{elsey_esedoglu_anisotropy, ejz}, with convolution kernels (i.e. the weights in weighted median filters) given by a variant of the kernel construction in \cite{ejz}.
The numerical studies presented suggest that when spatial resolution difficulties of threshold dynamics algorithms are alleviated (by considering their natural level set versions), there is good evidence of convergence to the desired dynamics using the approach to kernel construction presented in \cite{ejz} that \textcolor{black}{furnishes} a kernel for the desired anisotropic surface tension and anisotropic mobility pair of each interface in the network, at least in the simpler context of wetting / dewetting problems.
\bigskip

\noindent {\bf Acknowledgements:}
The authors were supported by NSF-DMS 2012015 and NSF-DMS 2410272.
\bibliographystyle{plain}
\bibliography{bib}
\end{document}